\documentclass[12pt,journal,12pt,draftclsnofoot,onecolumn]{IEEEtran}
\usepackage{tgtermes}
\usepackage{tgheros}
\usepackage{tgcursor}
\usepackage[T1]{fontenc}
\usepackage[utf8]{inputenc}
\usepackage{float}
\usepackage{amsthm}
\usepackage{amsmath}
\usepackage{amssymb}
\usepackage{graphicx}
\usepackage{setspace}
\providecommand{\tabularnewline}{\\}
\floatstyle{ruled}
\newfloat{algorithm}{tbp}{loa}
\providecommand{\algorithmname}{Algorithm}
\floatname{algorithm}{\protect\algorithmname}

\numberwithin{figure}{section}
\theoremstyle{plain}
\newtheorem{thm}{\protect\theoremname}
\newtheorem{prop}[thm]{\protect\propositionname}

\usepackage{float}
\usepackage{amsfonts}
\usepackage{amsthm}
\usepackage{amsmath}
\usepackage{amssymb}
\usepackage{graphicx}

\usepackage{tikz}
\usepackage{tikzscale}

\usepackage{gnuplot-lua-tikz}

\usepackage{fouriernc,bm} 

\usepackage{blindtext}   
\usepackage{adjustbox}

\usepackage{environ}

\usepackage{epsfig}
\usepackage{cite}

\usepackage{cancel}

\usepackage{epsfig}
\usepackage{mathrsfs}
\usepackage{algorithm}
\usepackage{algorithmic}
\usepackage{graphicx}

\usepackage{epsfig}

\usepackage{algorithmic}
\usepackage{color}

\usepackage{tocloft}

\usepackage{ragged2e}
\usepackage{listings}
\usepackage{xcolor}

\usepackage{nccparskip}
\usepackage{subfig}

\makeatletter
\providecommand{\tabularnewline}{\\}
\floatstyle{ruled}
\newfloat{algorithm}{tbp}{loa}
\providecommand{\algorithmname}{Algorithm}
\floatname{algorithm}{\protect\algorithmname}

\numberwithin{figure}{section}
\theoremstyle{plain}
\theoremstyle{plain}

\makeatother

\def\E{\mathbb{E}}

\setkeys{Gin}{width=\linewidth}

\allowdisplaybreaks

\numberwithin{figure}{section}

\newcommand{\listappendixname}{LIST OF APPENDICES}

\newcommand{\listofappendices}{%
\chapter*{\listappendixname}\@starttoc{app}}

\newcommand{\bigO}{\mathrm{O}}

\def\R{\mathbb{R}}

\def\H{\mathbb{H}}
\def\E{\mathbb{E}}

\allowdisplaybreaks

\def\R{\mathbb{R}}

\def\H{\mathbb{H}}
\def\E{\mathbb{E}}

\newcommand{\vct}[1]{\bm{#1}}

\SetParskip{1.1pt plus 1pt minus .4pt}
\allowdisplaybreaks 

\renewcommand{\contentsname}{\Large \center{TABLE OF CONTENTS}}

\newcommand{\norm}[1]{\left\lVert#1\right\rVert}

\providecommand{\propositionname}{Proposition}
\providecommand{\theoremname}{Theorem}

\begin{document}

\title{Deterministic Versus Randomized Kaczmarz Iterative Projection}

\author{\IEEEauthorblockN{Tim Wallace} and \IEEEauthorblockN{Ali Sekmen}\IEEEauthorrefmark{2}\\
\thanks{ \protect\IEEEauthorrefmark{2}\protect\IEEEauthorblockA{Email: asekmen@tnstate.edu}%
}Department of Computer Science\\
Tennessee State University\\
Nashville, TN USA}
\maketitle
\begin{abstract}
\begin{singlespace}
\noindent The Kaczmarz's alternating projection method has been widely
used for solving a consistent (mostly over-determined) linear system
of equations $A\bm{x}=\bm{b}$. Because of its simple iterative nature
with light computation, this method was successfully applied in computerized
tomography. Since tomography generates a matrix $A$ with highly coherent
rows, randomized Kaczmarz algorithm is expected to provide faster
convergence as it picks a row for each iteration at random, based
on a certain probability distribution. It was recently shown that
picking a row at random, proportional with its norm, makes the iteration
converge exponentially in expectation with a decay constant that depends
on the scaled condition number of $A$ and not the number of equations.
Since Kaczmarz's method is a subspace projection method, the convergence
rate for simple Kaczmarz algorithm was developed in terms of subspace
angles. This paper provides analyses of simple and randomized Kaczmarz
algorithms and explain the link between them. It also propose new
versions of randomization that may speed up convergence.\end{singlespace}

\end{abstract}

\section{Introduction}

Kaczmarz (in~\cite{Kaczmarz1937}) introduced an iterative algorithm
for solving a consistent linear system of equations $A\bm{x}=\bm{b}$
with $A\in\R^{M\times N}$. This method projects the estimate $\bm{x}^{j}$
onto a subspace normal to the row $a_{i}$ at step $j+1$ cyclically
with $i=j\pmod{M}+1$. The block Kaczmarz algorithm first groups the
rows into matrices $A_{1},A_{2},\hdots,A_{k}$ and then it projects
the estimate $\bm{x}^{j}$ onto the subspace normal to the subspace
spanned by the rows of $A_{i}$ at step $j+1$ cyclically with $i=j\pmod{k}+1$.
Obviously, the block Kaczmarz is equivalent to the simple Kaczmarz
for $k=M$. The Kaczmarz method is a method of alternating projection
(MAP) and it has been widely used in medical imaging as an algebraic
reconstruction technique (ART)~\cite{Gordon1970,Herman2009} due
to its simplicity and light computation. Strohmer \emph{et al.}~\cite{Strohmer2009}
proved that if a row for each iteration is picked in a random fashion
with probability proportional with $\ell_{2}$ norm of that row, then
the algorithm converges in expectation exponentially with a rate that
depends on a scaled condition number of $A$ (not on the number of
equations). Needell (in \cite{Needell2010}) extended the work of
\cite{Strohmer2009} for noisy linear systems and developed a bound
for convergence to the least square solution for $A\bm{x}=\bm{b}$.
Needell also developed a randomized Kaczmarz method that improves
the incoherency for iteration~\cite{Needell2013} and she analyzed
the convergence of randomized block Kaczmarz method~\cite{Needell2014}.
Chen and Powell (in~\cite{Chen2012}) consider a random measurement
matrix $A$ instead of random selection of measurements. Galantai
(in~\cite{Galantai2004,Galantai2005}) provides convergence analysis
for block Kaczmarz method by expanding the convergence analysis (based
on subspace angles) of Deutsch~\cite{Deutsch1997}. Brezinski (in~\cite{Brezinski2013})
utilizes the work of Galantai for accelerating convergence of regular
Kaczmarz method. 

\subsection{Paper Contributions}
\begin{itemize}
\item Research on regular and randomized Kaczmarz methods appear disconnected
in the literature. Even though convergence rates have been studied
separately, the link between them has not been explored sufficiently. 
\item A new randomization technique based on subspace angles has been developed
which indicates an advantage with coherent data measurements. 
\item A further method is introduced which orthogonalizes the subspace blocks
in order to mitigate the coherency. Convergence is consistent with
statistical expectations from theory and simulations.
\item The effects of measurement coherence are observed in the literature
and illustrated in our simulations with norm and angle based iteration
randomization.
\item A broader review and mathematical analysis of common methods is presented
from both statistical and deterministic perspectives. 
\end{itemize}

\section{Convergence of Regular Block Kaczmarz Method}

Let $\bm{x}^{*}$ be the solution of consistent $A\bm{x}=\bm{b}$
where $A\in\R^{M\times M}$ is full column rank. Let $A$ be row-partitioned
as $\{A_{1},\ldots,A_{k}\}$ where $A_{i}\in\R^{M_{i}\times M}$.
Then, the simple block Kaczmarz update is as follows: 
\begin{equation}
\bm{x}_{j+1}=\bm{x}_{j}+A_{i}^{T}(A_{i}A_{i}^{T})^{-1}(\bm{b}_{i}-A_{i}\bm{x}_{j})\;\;\; i=j\pmod{k}+1
\end{equation}
where $\bm{b}_{i}$ is the section of $\bm{b}$ that corresponds to
the rows of $A_{i}$. Note that since $A_{i}$ is full row rank, $A_{i}^{T}(A_{i}A_{i}^{T})^{-1}$
is the right pseudo-inverse of $A_{i}$. This is equivalent to: 
\begin{align*}
\bm{x}_{j+1} & =\bm{x}_{j}+A_{i}^{T}(A_{i}A_{i}^{T})^{-1}(A_{i}\bm{x}^{*}-A_{i}\bm{x}_{j})\\
\bm{x}_{j+1}-\bm{x}^{*} & =\bm{x}_{j}-\bm{x}^{*}-A_{i}^{T}(A_{i}A_{i}^{T})^{-1}A_{i}(\bm{x}_{j}-\bm{x}^{*}).
\end{align*}
Note that $A_{i}^{T}(A_{i}A_{i}^{T})^{-1}A_{i}$ is the projection
matrix for projection of the range of $A_{i}^{T}$: 
\begin{align}
\bm{x}_{j+1}-\bm{x}^{*} & =\bm{x}_{j}-\bm{x}^{*}-P_{Sp(A_{i}^{T})}(\bm{x}_{j}-\bm{x}^{*})\label{eqn:proj-1}\\
\bm{x}_{j+1}-\bm{x}^{*} & =(I-P_{Sp(A_{i}^{T})})(\bm{x}_{j}-\bm{x}^{*})\nonumber \\
\bm{x}_{j+1}-\bm{x}^{*} & =P_{Sp^{\perp}(A_{i}^{T})}(\bm{x}_{j}-\bm{x}^{*}).\label{eqn:projperp-1}
\end{align}

For one cycle of the blocks, 
\begin{align}
\bm{x}_{k}-\bm{x}^{*} & =P_{Sp^{\perp}(A_{k}^{T})}P_{Sp^{\perp}(A_{k-1}^{T})}\hdots P_{Sp^{\perp}(A_{1}^{T})}(\bm{x}_{0}-\bm{x}^{*}).\label{eqn:genproj-1}
\end{align}

Note that if $A\in\R^{M\times N}$ is a full column rank with $M<N$,
then the simple block Kaczmarz update is as follows: 
\begin{equation}
\bm{x}_{j+1}=\bm{x}_{j}+A_{i}^{\dagger}(\bm{b}_{i}-A_{i}\bm{x}_{j})=\bm{x}_{j}+A_{i}^{\dagger}A_{i}(\bm{x}^{*}-\bm{x}_{j})\;\;\; i=j\pmod{k}+1
\end{equation}
where $A_{i}^{\dagger}$ is the pseudo-inverse of $A_{i}$ and $A_{i}^{\dagger}A_{i}$
is the orthogonal projection onto $Sp(A_{i}^{T})$. Then, we get the
same equation as Equation \eqref{eqn:proj-1}, and subsequently we
get Equation \eqref{eqn:genproj-1}, 
\begin{align}
\bm{x}_{j+1}-\bm{x}^{*} & =\bm{x}_{j}-\bm{x}^{*}-P_{Sp(A_{i}^{T})}(\bm{x}_{j}-\bm{x}^{*}).
\end{align}

\subsection{Exponential Convergence}
\begin{thm}
Let $\bm{x}^{*}$ be the solution of consistent $A\bm{x}=\bm{b}$
where $A\in\R^{M\times M}$ is full column rank. Let $A$ be row-partitioned
as $\{A_{1},\ldots,A_{k}\}$ where $A_{i}\in\R^{M_{i}\times M}$.
Then, the simple block Kaczmarz converges exponentially and the convergence
rate depends of the number of blocks.\end{thm}
\begin{IEEEproof}
By Equation \eqref{eqn:proj-1} and orthogonal projection, 
\begin{equation}
\norm{\bm{x}_{j+1}-\bm{x}^{*}}_{2}^{2}=\norm{\bm{x}_{j}-\bm{x}^{*}}_{2}^{2}-\norm{P_{Sp(A_{i}^{T})}(\bm{x}_{j}-\bm{x}^{*})}_{2}^{2}.\label{eqn:diff-3}
\end{equation}
So, 
\begin{equation}
\norm{\bm{x}_{j+1}-\bm{x}^{*}}_{2}^{2}\leq\norm{\bm{x}_{j}-\bm{x}^{*}}_{2}^{2},\label{eqn:norm-1}
\end{equation}
$\bm{x}_{j}-\bm{x}^{*}$ depends on the initial condition $\tilde{\bm{x}}_{0}=\bm{x}_{0}-\bm{x}^{*}$,
and this dependence is scale-invariant. To see this, let $\bm{e}_{j}=\bm{x}_{j}-\bm{x}^{*}$
and consider $c\tilde{\bm{x}}_{0}$ where $c\in\R$. By Equation~\eqref{eqn:projperp-1},
\begin{align}
\bm{e}_{j+1}(c\tilde{\bm{x}}_{0}) & =P_{Sp^{\perp}(A_{j+1}^{T})}\bm{e}_{j}(c\tilde{\bm{x}}_{0})\nonumber \\
 & =P_{Sp^{\perp}(A_{j+1}^{T})}P_{Sp^{\perp}(A_{j}^{T})}\hdots P_{Sp^{\perp}(A_{1}^{T})}\bm{e}_{0}(c\tilde{\bm{x}}_{0})\nonumber \\
 & =P_{Sp^{\perp}(A_{j+1}^{T})}P_{Sp^{\perp}(A_{j}^{T})}\hdots P_{Sp^{\perp}(A_{1}^{T})}(c\tilde{\bm{x}}_{0})\nonumber \\
 & =cP_{Sp^{\perp}(A_{j+1}^{T})}P_{Sp^{\perp}(A_{j}^{T})}\hdots P_{Sp^{\perp}(A_{1}^{T})}\bm{e}_{0}(\tilde{\bm{x}}_{0})\nonumber \\
 & =c\bm{e}_{j+1}(\tilde{\bm{x}}_{0}).\label{eqn:scale-1}
\end{align}

We will first show that if $\bm{x}_{0}\neq\bm{x}^{*}$, then $\norm{\bm{x}_{k}-\bm{x}^{*}}_{2}<\norm{\bm{x}_{0}-\bm{x}^{*}}_{2}$.
By the way of contradiction, assume that $\bm{x}_{0}\neq\bm{x}^{*}$
and $\norm{\bm{x}_{k}-\bm{x}^{*}}_{2}=\norm{\bm{x}_{0}-\bm{x}^{*}}_{2}$.
By Equation~\eqref{eqn:norm-1}, 
\[
\norm{\bm{x}_{k}-\bm{x}^{*}}_{2}\leq\norm{\bm{x}_{k-1}-\bm{x}^{*}}_{2}\hdots<\norm{\bm{x}_{0}-\bm{x}^{*}}_{2}
\]
and therefore $\norm{\bm{x}_{l}-\bm{x}^{*}}_{2}=\norm{\bm{x}_{0}-\bm{x}^{*}}_{2}$
for all $1\leq l\leq k$. By Equation~\eqref{eqn:proj-1}, $P_{Sp(A_{l}^{T})}(\bm{x}_{l-1}-\bm{x}^{*})=0$
for all $1\leq l\leq k$. By Equation~\eqref{eqn:diff-3}, we get
$\bm{x}_{l}=\bm{x}_{0}$ for all $1\leq l\leq k$. This implies that
$P_{Sp(A_{l}^{T})}(\bm{x}_{0}-\bm{x}^{*})=0$ for all $1\leq l\leq k$.
So, 
\begin{align*}
P_{Sp^{\perp}(A_{k}^{T})\cap Sp^{\perp}(A_{k}^{T})\hdots\cap Sp^{\perp}(A_{1}^{T})}(\bm{x}_{0}-\bm{x}^{*}) & =0\\
P_{Sp^{\perp}(A^{T})}(\bm{x}_{0}-\bm{x}^{*}) & =0.
\end{align*}
Since $A$ is full column rank we get $\bm{x}_{0}=\bm{x}^{*}$, which
is a contradiction. So we know that $\norm{\bm{x}_{k}-\bm{x}^{*}}_{2}<\norm{\bm{x}_{0}-\bm{x}^{*}}_{2}$
(for one full cycle of $k$-iterations).

By compactness, there exists an $\epsilon\in(0,1)$ such that for
all $\tilde{\bm{x}}_{0}=\bm{x}_{0}-\bm{x}^{*}\in S^{N-1}$, 
\begin{equation}
\norm{\bm{x}_{k}-\bm{x}^{*}}_{2}\leq1-\epsilon.\label{eqn:epsilon-1}
\end{equation}
By Equations \eqref{eqn:scale-1} and \eqref{eqn:epsilon-1} 
\begin{align*}
\norm{\bm{x}_{k}-\bm{x}^{*}}_{2} & =\norm{\tilde{\bm{x}}_{0}}_{2}\bm{e}_{k}(\frac{\tilde{\bm{x}}_{0}}{\norm{\tilde{\bm{x}}_{0}}_{2}})\leq(1-\epsilon)\norm{\tilde{\bm{x}}_{0}}_{2}\\
\norm{\bm{x}_{k}-\bm{x}^{*}}_{2} & \leq(1-\epsilon)\norm{\bm{x}_{0}-\bm{x}^{*}}_{2}.
\end{align*}
Now consider iteration for $q$ cycles, 
\begin{align*}
\norm{\bm{x}_{qk}-\bm{x}^{*}}_{2} & \leq(1-\epsilon)^{q}\norm{\bm{x}_{0}-\bm{x}^{*}}_{2}\\
\norm{\bm{x}_{qk}-\bm{x}^{*}}_{2} & \leq[(1-\epsilon)^{1/k}]^{qk}\norm{\bm{x}_{0}-\bm{x}^{*}}_{2}.
\end{align*}
Therefore, we conclude that the exponential decay depends on the number
of blocks $k$. Note that $k=M$ for regular simple Kaczmarz and the
exponential decay depends on the number of rows in this case. The
randomized Kaczmarz algorithm proposed by Strohmer and Vershynin~\cite{Strohmer2009}
avoids this and it converges in expectation as $\E\norm{\bm{x}_{p}-\bm{x}^{*}}_{2}^{2}\leq(1-\kappa(A)^{-2})^{p}\norm{\bm{x}_{0}-\bm{x}^{*}}_{2}^{2}$,
where $\kappa(A)=\norm{A}_{F}\norm{A^{\dagger}}_{2}$ is the scaled
condition number of matrix $A$ with $A^{\dagger}$ is the pseudo-inverse
of $A$.
\end{IEEEproof}

\subsection{Iterative Subspace Projection Approach}

We can use the following theorem (in \cite{Galantai2005,Deutsch1997})
to show the convergence of regular block Kaczmarz method.
\begin{thm}
\label{Halperin}Let $M_{1},M_{2},\hdots M_{k}$ be closed subspaces
of the real Hilbert space $\H$. Let $M=\cap_{i=1}^{k}M_{i}$ and
$P_{M_{i}}\;(i=1,\hdots,k)$ be orthogonal projection on $M_{i}$.
Then, for each $\bm{x}\in\H$, 
\[
\lim_{q\to\infty}(P_{M_{k}}P_{M_{k-1}}\hdots P_{M_{1}})^{q}\bm{x}=P_{M}\bm{x}
\]
 where $P_{M}$ is the orthogonal intersection projection. 
\end{thm}
The block Kaczmarz is an alternating projection method with $M_{1}=Sp^{\perp}(A_{1}^{T}),\hdots,M_{k}=Sp^{\perp}(A_{k}^{T})$.
Also, $P_{M_{1}}=P_{Sp^{\perp}(A_{1}^{T})},\hdots,P_{M_{k}=Sp^{\perp}(A_{k}^{T})}$
and $M=Sp^{\perp}(A_{1}^{T})\cap\hdots\cap Sp^{\perp}(A_{k}^{T})=Sp^{\perp}(A^{T})$.
Since $A$ is full column rank, $Sp^{\perp}(A^{T})=\{0\}$ and $P_{M}=\{0\}$.
After $q$ cycles, 
\begin{equation}
\bm{x}_{qk}-\bm{x}^{*}=(P_{M_{k}}P_{M_{k-1}}\hdots P_{M_{1}})^{q}(\bm{x}_{0}-\bm{x}*).
\end{equation}
By Theorem~\ref{Halperin}, $\lim_{q\to\infty}\bm{x}_{qk}-\bm{x}^{*}=0$
and $\lim_{q\to\infty}\bm{x}_{qk}=\bm{x}^{*}$. Galantai in \cite{Galantai2005}
gives a bound for $\norm{\bm{x}_{qk}-\bm{x}^{*}}_{2}$ in terms of
principle angles between $M_{i}$\rq{}s.

\subsection{Bound for Block Kaczmarz in terms of Principle Angles}

Smith, Salmon, and Wagner established the following convergence theorem
for applying the alternating projection method in tomography \cite{Galantai2005,smith1977}: 
\begin{thm}
\label{angles} Let $M_{1},M_{2},\hdots M_{k}$ be closed subspaces
of the real Hilbert space $\H$. Let $M=\cap_{i=1}^{k}M_{i}$ and
$P_{M_{i}}\;(i=1,\hdots,k)$ be orthogonal projection on $M_{i}$
($P_{M}$ is the orthogonal intersection projection). Let $\theta_{j}=\alpha(M_{j},\cap_{i=j+1}^{k}M_{i})$,
then for each $\bm{x}\in\H$ and integer $q\geq1$, 
\[
\norm{(P_{M_{k}}P_{M_{k-1}}\hdots P_{M_{1}})^{q}\bm{x}-P_{M}\bm{x}}_{2}^{2}\leq(1-\Pi_{j=1}^{k-1}\sin^{2}\theta_{j})^{q}\norm{\bm{x}-P_{M}\bm{x}}_{2}^{2}
\]
where $P_{M}$ is the orthogonal intersection projection. 
\end{thm}
In the special case of the block Kaczmarz, we have $\H=\R^{N}$, $M_{1}=Sp^{\perp}(A_{1}^{T}),\hdots,M_{k}=Sp^{\perp}(A_{k}^{T})$.
Also, $P_{M_{1}}=P_{Sp^{\perp}(A_{1}^{T})},\hdots,P_{M_{k}}=P_{Sp^{\perp}(A_{k}^{T})}$
and $M=Sp^{\perp}(A_{1}^{T})\cap\hdots\cap Sp^{\perp}(A_{k}^{T})=Sp^{\perp}(A^{T})$.
Since $A$ is full column rank, $Sp^{\perp}(A^{T})=\{0\}$ and $P_{M}=\{0\}$.
Therefore, after $q$ cycles, 
\begin{equation}
\norm{\bm{x}_{qk}-\bm{x}^{*}}_{2}^{2}=\norm{(P_{M_{k}}P_{M_{k-1}}\hdots P_{M_{1}})^{q}(\bm{x}_{0}-\bm{x}^{*})}_{2}^{2}\leq(1-\Pi_{j=1}^{k-1}\sin^{2}\theta_{j})^{q}\norm{\bm{x}_{o}-\bm{x}^{*}}_{2}^{2}
\end{equation}
where $\theta_{j}$ is as defined in Theorem~\ref{angles}. Note
that the exponential decay rate depends on the number of blocks $k$
as shown below. 
\begin{equation}
\norm{\bm{x}_{qk}-\bm{x}^{*}}_{2}^{2}\leq[(1-\Pi_{j=1}^{k-1}\sin^{2}\theta_{j})^{1/k}]^{qk}\norm{\bm{x}_{o}-\bm{x}^{*}}_{2}^{2}
\end{equation}
Galantai in~\cite{Galantai2005} developed another bound (for $A\in\R^{M\times M}$)
by defining a new matrix $X_{i}$ for each block $A_{i}$ as follows:
\begin{thm}
\label{angles2}Let $\bm{x}^{*}$ be the solution of $A\bm{x}=\bm{b}$
for a consistent linear system with $A\in\R^{M\times M}$. Let $A$
be row-partitioned as $\{A_{1},\ldots,A_{k}\}$ where $A_{i}\in\R^{M_{i}\times N}$.
Let $M_{1}=Sp^{\perp}(A_{1}^{T}),\hdots,M_{k}=Sp^{\perp}(A_{k}^{T})$
and $A_{i}A_{i}^{T}=LL^{T}$ be the Cholesky decomposition of $A_{i}A_{i}^{T}$.
Define $X_{i}=A_{i}^{T}L^{-T}$ and $X=[X_{1},\hdots,X_{k}]$. Then
for each $\bm{x}\in\R^{N}$ and integer $q\geq1$, 
\[
\norm{\bm{x}_{qk}-\bm{x}^{*}}_{2}^{2}\leq[1-\det(X^{T}X)]^{q}\norm{\bm{x}_{o}-\bm{x}^{*}}_{2}^{2}=[(1-\det(X^{T}X))^{1/k}]^{qk}\norm{\bm{x}_{o}-\bm{x}^{*}}_{2}^{2}
\]

\end{thm}

\subsection{Special Case: Simple Kaczmarz for $A\in\R^{M\times M}$}

Note that this section assumes that $A\in\R^{M\times M}$. The block
Kaczmarz algorithm is equivalent to the simple Kaczmarz algorithm
if the number of blocks $k$ is equal to the number of rows $M$.
In this case, $A_{i}A_{i}^{T}=\norm{\bm{a}_{i}}_{2}^{2}=LL^{T}$.
therefore, $L=\norm{\bm{a}_{i}}_{2}$ and $L^{-T}=1/\norm{\bm{a}_{i}}_{2}$.
This implies that $X_{i}=[\frac{\bm{a}_{i}}{\norm{\bm{a}_{i}}_{2}}]$.
Then, $X\in\R^{M\times M}$ is defined as: 
\begin{equation}
X=[\frac{\bm{a}_{1}}{\norm{\bm{a}_{1}}_{2}},\hdots,\frac{\bm{a}_{M}}{\norm{\bm{a}_{M}}_{2}}].
\end{equation}
Assume the matrix $A$ has normalized rows and we pick a row at each
iteration uniformly randomly. Note that this assumption is feasible
as scaling a row of $A$ and the corresponding measurement in $\bm{b}$
does not change the solution $\bm{x}$.

$X$ is the Gram matrix with $0\leq\det(X^{T}X)\leq\norm{\bm{x}_{1}}_{2}^{2}\norm{\bm{x}_{2}}_{2}^{2}\hdots\norm{\bm{x}_{M}}_{2}^{2}$.
Since $\norm{\bm{x}_{i}}_{2}=1$ and $X$ is full rank, we have $0<\det(X^{T}X)\leq1$.
Using Theorem~\ref{angles2}, we get the following deterministic
bound: 
\begin{equation}
\norm{\bm{x}_{qM}-\bm{x}^{*}}_{2}^{2}\leq[(1-\det(X^{T}X))^{1/M}]^{qM}\norm{\bm{x}_{0}-\bm{x}^{*}}_{2}^{2}.
\end{equation}
Since $A$ is normalized, we get, $X=A^{T}$ and therefore: 
\begin{equation}
\norm{\bm{x}_{qM}-\bm{x}^{*}}_{2}^{2}\leq[(1-\det(AA^{T}))^{1/M}]^{qM}\norm{\bm{x}_{0}-\bm{x}^{*}}_{2}^{2}.\label{eqn:normal}
\end{equation}
Bai \emph{et al.} (in \cite{Bai2013}) uses the Meany Inequality to
develop a general form of this inequality.

\section{Randomized Kaczmarz Method}

\subsection{Randomization Based on Row $\ell_{2}$ Norms}

Strohmer \emph{et al.} (in~\cite{Strohmer2009}) developed a randomized
Kaczmarz algorithm that picks a row of $A$ in a random fashion with
probability proportional with $\ell_{2}$ norm of that row. They proved
that this method has exponential expected convergence rate. Since
the rows are picked based on a probability distribution generated
by the $\ell_{2}$ norms of the rows of $A$, it is clear that scaling
some of the equations does not change the solution set. However, it
may drastically change the order of the rows picked at each iteration.
Censor \emph{et al.} discusses (in~\cite{Censor09}) that this should
not be better than the simple Kaczmarz as picking a row based on its
$\ell_{2}$ norm does not change the geometry of the problem. Theorem
\ref{thm:randomkaczmarz} is from \cite{Strohmer2009}.
\begin{algorithm}
\protect\caption{Randomized Kaczmarz (of \cite{Strohmer2009})}

\begin{algorithmic}[1]\REQUIRE\label{2StepRK-1} \label{alg:random}An
over-determined linear set of consistent equations $A\bm{x}=\bm{b}$,
where $A$ is $M\times N$ matrix and $\bm{b}\in\R^{M}$. Let $\bm{a}_{1},\hdots,\bm{a}_{M}$
be the rows of $A$ and $b_{j}$ be the $j^{th}$ element of $\bm{b}$.
\STATE Pick an arbitrary initial approximation $\bm{x}_{0}$. \STATE
Set $p=0$. \WHILE {not converged} \STATE Randomly choose $r(i)$
from $\left\{ 1,\hdots,M\right\} $ with probability proportional
to $\norm{\bm{a}_{r(i)}}_{2}^{2}$. \STATE $\bm{x}_{p+1}=\bm{x}_{p}+\frac{b_{r(i)}-\langle\bm{a}_{r(i)},\bm{x}_{p}\rangle}{\norm{\bm{a}_{r(i)}}_{2}^{2}}\bm{a}_{r(i)}$
\STATE Set $p=p+1$ \ENDWHILE \end{algorithmic} 
\end{algorithm}

\begin{thm}
\label{thm:randomkaczmarz} Let $\bm{x}^{*}$ be the solution of $A\bm{x}=\bm{b}$
Then, Algorithm~\ref{alg:random} converges to $\bm{x}^{*}$ in expectation,
with the average error 
\begin{equation}
\E\norm{\bm{x}_{p}-\bm{x}^{*}}_{2}^{2}\leq(1-\kappa(A)^{-2})^{p}\norm{\bm{x}_{0}-\bm{x}^{*}}_{2}^{2}
\end{equation}
where $\kappa(A)=\norm{A}_{F}\norm{A^{\dagger}}_{2}$ is the scaled
condition number of matrix $A$ with $A^{\dagger}$ is the left pseudo-inverse
of $A$. 
\end{thm}
Note that $A$ is a full column matrix ($A\in\R^{M\times N}$ with
$rank(A)=N$) and therefore we define $A^{\dagger}$ as left pseudo-inverse
of $A$. We observe that the randomization should work better than
the simple (cyclic) Kaczmarz algorithm for matrices with highly coherent
rows (e.g. matrices generated by the computerized tomography). Since
the Kaczmarz algorithm is based on projections, the convergence will
be slow if the consecutive rows selected are highly coherent (i.e.
the angle between $\bm{a}_{i}$ and $\bm{a}_{i+1}$ is small). Picking
rows randomly (not necessarily based on the $\ell_{2}$ norms) makes
picking more incoherent rows possible in each iteration. Therefore,
the randomization may be useful for certain applications such as medical
imaging. Note that matrix $A$ generated by computerized tomography
has coherent and sparse rows due to physical nature of data collection.
In fact, using Theorem~\ref{thm:randomkaczmarz}, we can develop
the following proposition. 
\begin{prop}
Let $A\bm{x}=\bm{b}$ be a consistent linear system of equations ($A\in\R^{M\times N}$)
and let $\bm{x}_{0}$ be an arbitrary initial approximation to the
solution of $A\bm{x}=\bm{b}$. For $k=1,2,\hdots$ compute 
\begin{equation}
\bm{x}_{p+1}=\bm{x}_{p}+\frac{b_{r(i)}-\langle\bm{a}_{r(i)},\bm{x}_{p}\rangle}{\norm{\bm{a}_{r(i)}}_{2}^{2}}\bm{a}_{r(i)}
\end{equation}
where $r(i)$ is chosen from the set $\{1,2,\hdots,M\}$ at random,
with \textbf{any probability distribution}. Let $\bm{x}^{*}$ be the
solution of $A\bm{x}=\bm{b}$. Then, 
\begin{equation}
\E\norm{\bm{x}_{p}-\bm{x}^{*}}_{2}^{2}\leq(1-\kappa(B)^{-2})^{p}\norm{\bm{x}_{0}-\bm{x}^{*}}_{2}^{2}
\end{equation}
where $\kappa(B)=\norm{B}_{F}\norm{B^{\dagger}}_{2}$ is the scaled
condition number of a matrix $B$ that is obtained by some row-scaling
of $A$. \end{prop}
\begin{IEEEproof}
This is due to the fact that, row-scaling of $A$ (with scaling of
the corresponding $b$) does not change the geometry of the problem
and we can scale the rows to generate any probability distribution.
In other words, we can obtain another matrix $B$ from $A$ by scaling
its rows in such a way that picking the rows of $B$ based on the
$\ell_{2}$ norms of the rows will be equivalent to picking the rows
of $A$ based on the chosen probability distribution. Therefore, clearly,
any randomization of the row selection will have exponential convergence,
however, the rate will depend on the condition number of another matrix.
For example, if we use uniform distribution, we can then normalize
each row to have matrix $B$ as follows and then pick the rows at
random with probability proportional to the norms of the rows. 
\begin{equation}
B=[\frac{\bm{a}_{1}}{\norm{\bm{a}_{1}}_{2}},\hdots,\frac{\bm{a}_{M}}{\norm{\bm{a}_{M}}_{2}}]^{T}.
\end{equation}

\end{IEEEproof}

\subsection{Randomization based on Subspace Angles}

Our approach iterates through the rows of $A$ based on a probability
distribution using the hyperplane (subspace) angles. Therefore, it
is immune to scaling or normalization. This approach first generates
a probability distribution based on the angles between the hyperplanes
(represented by the rows of $A\bm{x}=\bm{b}$). Then, it randomly
picks two hyperplanes using this probability distribution. This is
followed by a two-step projection on these hyperplanes (see Algorithm~\ref{2StepRK}).
\begin{algorithm}
\protect\caption{Randomized Kaczmarz Hyperplane Angles}

\label{2StepRK} \begin{algorithmic}[1] \REQUIRE An over-determined
linear set of consistent equations $A\bm{x}=\bm{b}$, where $A$ is
$M\times N$ matrix and $\bm{b}\in\R^{M}$. Let $\bm{a}_{1},\hdots,\bm{a}_{M}$
be the rows of $A$ and $\bm{b}_{j}$ be the $j^{th}$ element of
$\bm{b}$. \STATE Pick an arbitrary initial approximation $\bm{x}_{0}$.
\STATE Set $k=0$. \STATE Randomly choose $f(i)$ from $\{1,2,\hdots,M\}$
with a uniform distribution. \WHILE {not converged} \STATE Randomly
choose $g(i)$ from $\left\{ 1,\hdots,M\right\} $ with probability
proportional to $1-\dfrac{\langle\bm{a}_{f(i)},\bm{a}_{g(i)}\rangle^{2}}{\|\bm{a}_{f(i)}\|_{2}^{2}\|\bm{a}_{g(i)}\|_{2}^{2}}$
\STATE Compute $\bm{x}_{k+1}=\bm{x}_{k}+\dfrac{b_{f(i)}-\langle\bm{a}_{f(i)},\bm{x}_{k}\rangle}{\|\bm{a}_{f(i)}\|_{2}^{2}}a_{f(i)}$
\STATE Compute $\bm{x}_{k+2}=\bm{x}_{k+1}+\dfrac{b_{g(i)}-\langle\bm{a}_{g(i)},\bm{x}_{k}\rangle}{\|\bm{a}_{g(i)}\|_{2}^{2}}a_{g(i)}$
\STATE Set $f(i)=g(i)$ \STATE Set $k=k+2$ \ENDWHILE \end{algorithmic} 
\end{algorithm}

\subsection{P-Subspaces Approach}

A new method has been developed which is intended to better accommodate
the coherency of non-orthogonal data measurements. This next section
makes contributions towards proving the statistical convergence of
the randomized Kaczmarz orthogonal subspace (RKOS) algorithm. As described
in \cite{Wallace2013}, the RKOS initially uses $\ell^{2}$-norm random
hyperplane selection and subsequent projection into a constructed
$P-$dimensional orthogonal subspace $S_{P}$ comprised of an additional
$P-1$ hyperplanes selected uniformly at random.

The algorithm uses a recursive method to solve for the projections
into the orthogonal subspace which is constructed using Gram-Schmidt
(GS) procedure. However, a second approach demonstrates an alternate
method of arriving at similar results, based upon an a closed form
matrix for QR decomposition \cite{golub-vanloan:1996} of projection
blocks.

In each of the above cases, vector operations inside the orthogonal
subspace preserve the $\ell^{2}$-norm, and reduce errors that would
normally be induced for coherent non-orthogonal projections which
may be present in the simple Kaczmarz.

\subsubsection{Orthogonal Subspaces}

A statistical convergence analysis for Randomized Kaczmarz Orthogonal
Subspace (RKOS) method is developed assuming identically and independently
distributed (IID) random variables as vector components of each row
of the measurement matrix $A$.

\paragraph{Orthogonal Construction }

In many problems, $M{\gg}N$ and fast but optimal solutions are needed,
often in noisy environments. In most cases, orthogonal data projection
sampling is not feasible due to the constraints of the measurement
system. The algorithm and procedure for the RKOS method is given in
reference \cite{Wallace2013} and is intended to construct orthogonal
measurements subspaces (see Algorithm \eqref{alg:Subspace-Kaczmarz-Projections-1-1}). 

The general technique is to solve using a constructed orthogonal basis
from a full rank set of linearly independent measurements in for each
subspace in Gram-Schmidt fashion \cite{Yanai:1414711,Meyer:2000:MAA:343374}. 

{\footnotesize{}}The subspace estimation may be computed as $P-$dimensional
subspace projection into the subspace orthonormal vector basis:\vspace{-5pt}
\begin{equation}
\bm{x}_{S_{P}}=\sum_{l=1}^{P}\langle\hat{\bm{u}}_{l},\bm{x}\rangle\hat{\bm{u}}_{l}.
\end{equation}
where $\bm{x}_{S_{P}}$ in $S_{P}\subseteq S_{N}$ subspace is the
$P-$dimensional solution approximation which becomes exact for $S_{P=N}$
for $\bm{x}_{S_{P=N}}\in\mathbb{R}^{N}$ in the noiseless, self-consistent,
case.%
\footnote{The $u$ vector with the hat symbol $\hat{u}$ indicates unit $\ell^{2}$-norm%
}

\paragraph{Modified Kaczmarz}

The standard Kaczmarz equation is essentially iterative projections
into a single subspace of dimension one; based upon the sampling hyperplanes,
these projections are often oblique, especially in highly-coherent
sampling.

The approach herein is motivated towards constructing an iterative
algorithm based upon Kaczmarz which may be accelerated while controlling
the potential projection errors and incurring reasonable computational
penalty. The algorithm is simply to add subspaces of larger dimensions.
Let
\begin{eqnarray}
\bm{x}-\bm{x}_{k+1} & = & \bm{x}-\bm{x}_{k}-\sum_{l=1}^{P}\langle\hat{\bm{u}}_{l},\bm{x}-\bm{x}_{k}\rangle\hat{\bm{u}}_{l}.
\end{eqnarray}
It is convenient to make a substitution as follows:
\begin{equation}
\vct{\bm{z}}_{k+1}=\bm{x}-\bm{x}_{k+1}.
\end{equation}
Using above substitution and orthonormal condition%
\footnote{It is worthwhile to note that in the problem setup, a fixed vector
is projected into a randomized $P$-dimensional subspace, where algebraic
orthogonality was used to obtain Equation \eqref{eq:zl2-1}. In the
this statistical treatment of the same equation, the expectation of
two random unit vectors vanishes for independent uncorrelated zero
mean probability distribution functions, providing the statistical
orthogonality on average satisfying \eqref{eq:zl2-1}.%
} $\langle\hat{\bm{u}}_{j},\hat{\bm{u}}_{k}\rangle=\delta_{j,k}$,
where the Kronecker $\delta_{j,k}=\begin{cases}
0 & \text{if }j\neq k\\
1 & \text{if }j=k,
\end{cases}$, find the $\ell^{2}$-norm squared of $\vct{\bm{z}}_{k+1}$:
\begin{equation}
\Vert\vct{\bm{z}}_{k+1}\Vert_{2}^{2}=\Vert\vct{\bm{z}}_{k}\Vert_{2}^{2}-\sum_{l=1}^{P}|\langle\hat{\bm{u}}_{l},\vct{\bm{z}}_{k}\rangle|^{2}.\label{eq:zl2-1}
\end{equation}
\begin{figure}
\centering{}\subfloat[{\small{}\label{Sta:Unit-sphere-2-D-1}IID Gaussian Unit Vector Image}]{\centering{}\includegraphics[width=0.28\columnwidth]{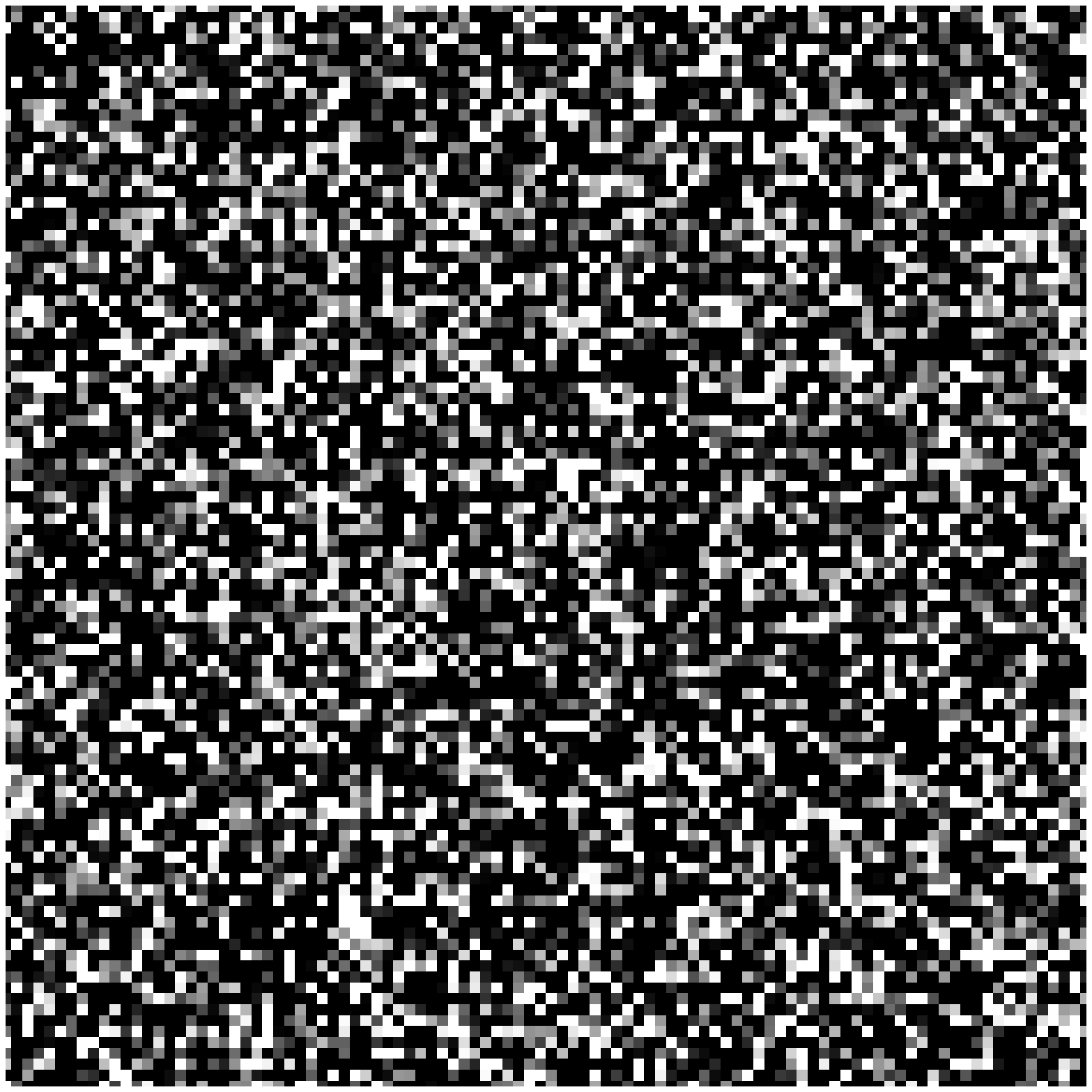}}
\ \ \ \ \subfloat[{\small{}\label{Sta:Phantom-image-1}CT Phantom Image}]{\centering{}\includegraphics[width=0.28\columnwidth]{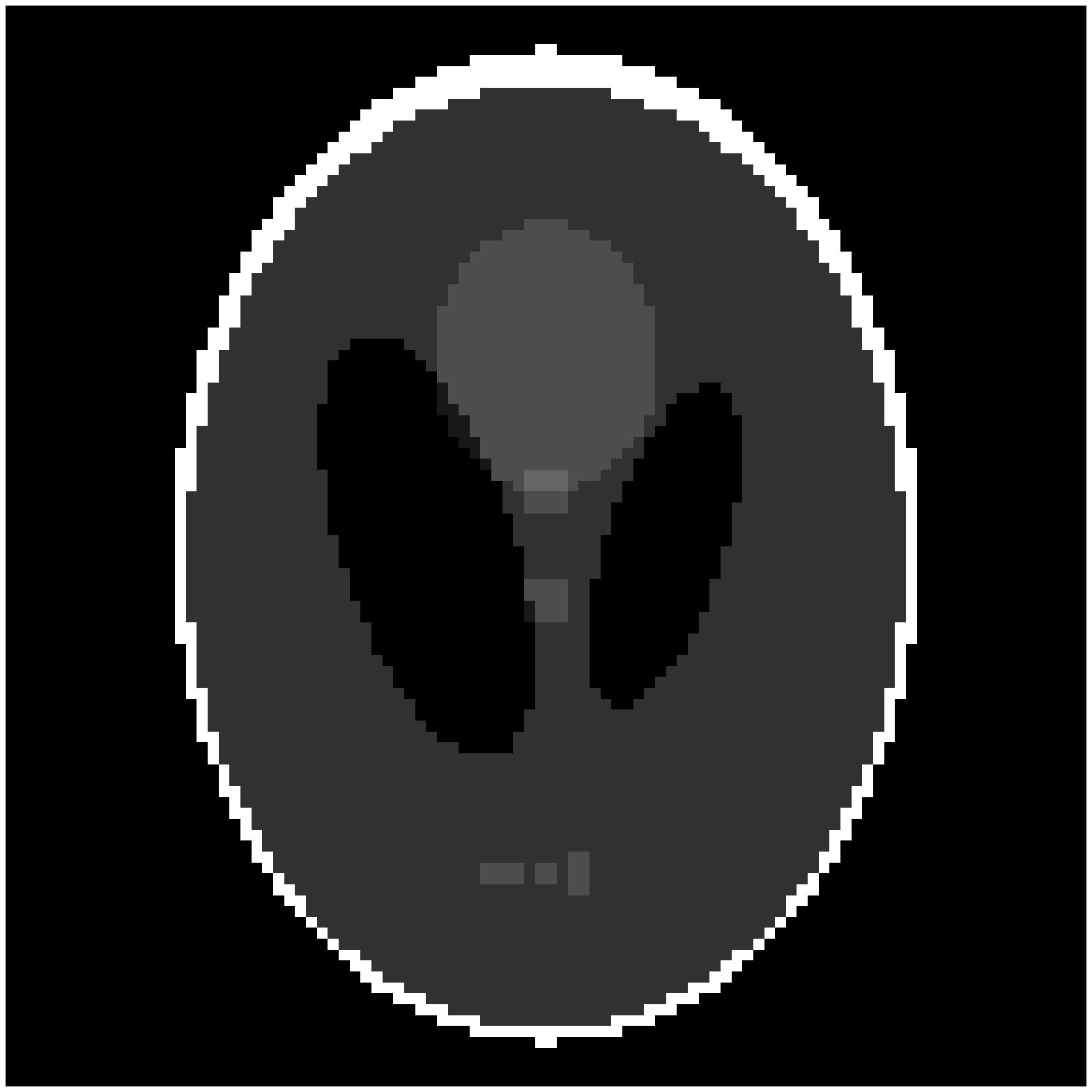}}\protect\caption{{\small{}Representative test data} }
\end{figure}
The ensemble average of the above Equation \ref{eq:zl2-1} yields
the convergence result, which is the main topic of this section. 
\begin{algorithm}
\protect\caption{\label{alg:Subspace-Kaczmarz-Projections-1-1}P-Subspace Kaczmarz
Projections}

\begin{algorithmic}[1]

\REQUIRE{\footnotesize{} }{\footnotesize \par}

{\footnotesize{}\vspace{4pt}
}{\footnotesize \par}

{Matrix $A\in\mathbb{R}^{M\times N}$ full-rank consistent measurements
subject to $A\bm{x}=\bm{b}$, for $\bm{b}\in\mathbb{R}^{M}$. }

\vspace{8pt}
\hrule\vspace{8pt}

\STATE Set $\bm{x}_{0}$ to initial approximation, $i=1$ \WHILE
{not converged} \STATE Select $\dim(S_{P})=P<N$ distinct linearly
independent rows of $A$ relative to random rule. Construct block
matrix $A_{i}\in\mathbb{R}^{P\times N}$ comprised of rows $\left\{ \vct{a}_{i,1},\ldots,\vct{a}_{i,P}\right\} $.\STATE
Perform Gram-Schmidt procedure on $A_{i}$ to obtain the orthonormal
set of columns $\left\{ \bm{u}_{i,1},\ldots,\bm{u}_{i,P}\right\} $.
Let $Q_{i}=\left\{ \bm{u}_{i,1},\ldots,\bm{u}_{i,P}\right\} \in\mathbb{R}^{N\times P}$
\STATE Update $\bm{x}_{i}$ as follows:

$\bm{x}_{i}=\bm{x}_{i-1}+Proj_{Sp(Q_{i})}(\bm{x}_{i-1}-\bm{x})$,

$\bm{x}_{i}=\bm{x}_{i-1}-Q_{i}Q_{i}^{T}(\bm{x}-\bm{x}_{i-1})$,

\STATE Compute $Q_{i}^{T}\bm{x}$ iteratively using $\left\{ \vct{a}_{i,1},\ldots,\vct{a}_{i,P}\right\} $,
$\left\{ \vct{b}_{i,1},\ldots,\vct{b}_{i,P}\right\} $, $\left\{ \bm{u}_{i,1},\ldots,\bm{u}_{i,P}\right\} $

\STATE Update $i=i+1$

\ENDWHILE \end{algorithmic} 
\end{algorithm}

\subsubsection{Convergence for IID Measurement Matrix}

Firstly, the expectation of a single random projection is computed.
In the second step, the terms are summed for the P-dimensional subspace.
Experimental results are included in a latter section.

\paragraph{Expectation of IID Projections}

\label{sec:Consider-the-expectationofIID}Consider the expectation
of the $\ell^{2}$-norm squared of the projection of fixed vector
$\bm{x}\in\mathbb{R}^{N\times1}$ onto a random subspace basis $U_{P}\in$
of dimension $P$, 
\[
\mathbb{E}[\|U_{P}^{T}\bm{x}\|_{2}^{2}],
\]
where the matrix basis $U_{P}\in\mathbb{R}^{NxP}$ is comprised of
$P-$columns of unit vectors $\hat{\bm{u}}_{j}\;\in\mathbb{R}^{N}$
in a constructed orthogonal basis for 
\begin{eqnarray}
\hat{\bm{u}}_{j} & \rightarrow & \hat{\bm{U}}_{j}=[U_{j,1},\ldots,U_{j,N}]\frac{1}{C_{\sigma}},\label{eq:unit random vector-1}\\
 & = & \frac{\bm{U}_{j}}{\|\bm{U}_{j}\|_{2}^{2}}\qquad\forall\, j\in\left[1,\ldots,P\right].
\end{eqnarray}
where the upper case components0 $U_{j,i}$ represent the $(j,i)$-th
IID random variable component, and normalization constant $C_{\sigma}$
is to be determined.

Further noting that complex conjugate $(.)^{*}$ reduces to transpose
$(.)^{T}$ for real components, the $\ell^{2}$-norm squared of the
projection expands to
\[
\|U_{P}^{T}\bm{x}\|_{2}^{2}=\bm{x}^{T}U_{P}U_{P}^{T}\bm{x}.
\]
In the next section, the goal is to find the expected value for outer
product of the projection, 
\[
\mathbb{E}\left[\bm{x}^{T}\hat{\bm{U}}_{j}\hat{\bm{U}}_{j}^{T}\bm{x}\right]\,\forall j\in[1,\ldots,P].
\]

\paragraph{Unit Vector }

The deterministic identity for the magnitude of a unit vector is well
known result for $\;\hat{\bm{u}}\in\mathbb{R}^{N}$, 
\begin{equation}
\|\hat{\bm{u}}\|_{2}^{2}=\sum_{i=1}^{N}\frac{u_{i}^{2}}{\|\bm{u}\|_{2}^{2}}=1.\label{eq:unit vector-2}
\end{equation}
The following statistical result must apply for the $j$-th column
unit vector:
\begin{eqnarray}
\mathbb{E}\left[\|\hat{\bm{U}}_{j}\|_{2}^{2}\right] & =\mathbb{E} & \left[\hat{\bm{U}}_{j}^{T}\hat{\bm{U}}_{j}\right]=1\label{eq:unit rand vec-1}\\
 & =\mathbb{E} & [U_{j,1}^{2}+\cdots+U_{j,N}^{2}]\frac{1}{C_{\sigma}^{2}}.\nonumber 
\end{eqnarray}

\paragraph{Normalization of Random Unit Vector }

Denote $\hat{\bm{U}}_{j}$ as the $j$-th random variable unit-norm
vector associated with a set of column vectors $\left\{ \bm{U}_{j}\right\} _{j\in1,\ldots,P}$
comprising a random subspace matrix $U_{N\times P}$ having IID random
variable components $U_{j,i}$. However, no additional assumptions
on the distribution of the random variables are made at this time,
other than IID. 

The expectation of both sides of Equation \eqref{eq:unit rand vec-1}
for random vector $\bm{U}_{j}$ are found such that:
\begin{equation}
\mathbb{E}\sum_{i=1}^{N}\frac{U_{j,i}^{2}}{C_{\sigma}^{2}}=\sum_{i=1}^{N}\mathbb{E}\left[\frac{U_{j,i}^{2}}{C_{\sigma}^{2}}\right]=1,
\end{equation}
\[
N\mathbb{\times}\frac{\mathbb{E}\left[U_{j,i}^{2}\right]}{C_{\sigma}^{2}}=1.
\]
Solving above for each unit vector component in this treatment implies
a random variable $U_{j,i}$ with zero mean and variance as follows:
\begin{equation}
\mathbb{E}\left[U_{j,i}^{2}\right]=\sigma_{j,i}^{2}=\frac{C_{\sigma}^{2}}{N}\qquad\forall\; U_{j,i\in1,\ldots,N}\in f(U_{j,i}),\label{eq:single-component-expectation-1}
\end{equation}
 where $f(U_{i,j})$ is the associated IID probability distribution.

\paragraph{P-Dimensional Random Projection}

The next step is to compute the expectation of the magnitude of the
projection of fixed vector $\bm{x}$ onto random $P$-dimensional
orthonormal subspace $U_{P}$ projection term by term. Let $\bm{\alpha}\in\mathbb{R}^{P}$
be a column vector defined as $\bm{\alpha}=U_{P}^{T}\bm{x}$ and find
the $\ell^{2}$-norm squared:
\begin{eqnarray}
\|\bm{\alpha}\|_{2}^{2} & = & \alpha_{1}^{2}+\alpha_{2}^{2}+\ldots+\alpha_{P}^{2}\\
 & = & \|U_{P}^{T}\bm{x}\|_{2}^{2}=\bm{x}^{T}U_{P}U_{P}^{T}\bm{x},\nonumber 
\end{eqnarray}
 where
\begin{eqnarray}
\alpha_{j}^{2} & = & \langle\hat{\bm{u}}_{j},\bm{x}\rangle^{2}\\
 & = & \left(u_{j,1}\bm{x}_{1}+\ldots+u_{j,N}\bm{x}_{N}\right)^{2}\\
 & = & \sum_{i,k}^{N,N}\frac{u_{j,k}u_{j,i}\bm{x}_{k}\bm{x}_{i}}{\|\bm{\bm{u}_{j}}\|_{2}^{2}}.
\end{eqnarray}

Let upper case $U_{j,k}$ denote the $k$-th IID element random%
\footnote{This is not the same k-variable as the Kaczmarz iteration variable%
} variable of the $j$-th column vector $\bm{U}_{j}$ associated with
column vector $\bm{u_{j}}$; let $\bm{x}$ vector denote a fixed point.
Next, take the expectation of the term over the possible outcomes
of $U_{j,k}$ random variables. Using the IID assumption, the expected
value for a single projection component preserves terms squared as
follows:
\begin{equation}
\mathbb{E}\left[\alpha_{j}^{2}\right]=\mathbb{E}\left[\sum_{i,k}^{N,N}\frac{U_{j,k}U_{j,i}\bm{x}_{k}\bm{x}_{i}}{C_{\sigma}^{2}}\right]=\sum_{i,k}^{N,N}\mathbb{E}\left[\frac{U_{j,k}U_{j,i}\bm{x}_{k}\bm{x}_{i}}{C_{\sigma}^{2}}\right]
\end{equation}
\[
=\sum_{k=1}^{N}\mathbb{E}\left[\frac{U_{k}^{2}\bm{x}_{k}^{2}}{C_{\sigma}^{2}}\right]=\sum_{k}^{N}\mathbb{E}\left[\frac{U_{j,k}^{2}}{C_{\sigma}^{2}}\right]\bm{x}_{k}^{2}
\]
\begin{eqnarray*}
 & = & \mathbb{E}\left[\frac{U_{j,k}^{2}}{C_{\sigma}^{2}}\right]\sum_{k}^{N}\bm{x}_{k}^{2}=\mathbb{E}\left[\frac{U_{j,k}^{2}}{C_{\sigma}^{2}}\right]\|\bm{x}\|_{2}^{2}\\
 & = & \frac{1}{C_{\sigma}^{2}}\frac{C_{\sigma}^{2}}{N}\|\bm{x}\|_{2}^{2}\\
 & = & \frac{1}{N}\|\bm{x}\|_{2}^{2}.
\end{eqnarray*}

It is now possible to determine the expectation for $P$-terms of
the projection as, 
\begin{equation}
\mathbb{E}\left[\|\bm{\alpha}\|_{2}^{2}\right]=\mathbb{E}\left[\sum_{j=1}^{P}\alpha_{j}^{2}\right]=\frac{P}{N}\|\bm{x}\|_{2}^{2}\label{eq:P/N E(|z|^2)-1}
\end{equation}
 subject to IID constraint on $\hat{\bm{U}}_{j}$ where it is further
noted that $\sigma^{2}N=C_{\sigma}^{2}$ in Equation \eqref{eq:single-component-expectation-1}.

\paragraph{Error per Iteration}

For a given $k$-th Kaczmarz iteration, the expectation of the projection
of fixed vector $\bm{x}$ onto the random P-dimensional subspace $U_{P}$
is known from above. The total convergence expectation may then be
computed, using a method similar to Strohmer's, starting%
\footnote{Recall that derivation of this equation \eqref{eq:z^2-1} requires
orthogonality among the $\hat{u}_{l}$ subspace basis vectors.%
} with Equation \eqref{eq:z^2-1}:
\begin{equation}
\|\bm{z}_{k+1}\|_{2}^{2}=\|\bm{z}_{k}\|_{2}^{2}-\sum_{l=1}^{P}|\langle\bm{z}_{k},\hat{\bm{u}}_{l}\rangle|^{2}\label{eq:z^2-1}
\end{equation}
\begin{equation}
\mathbb{E}_{\left\{ k+1|\bm{z}_{0},\bm{z}_{1},\ldots,\bm{z}_{k}\right\} }\left[\|\bm{z}_{k+1}\|_{2}^{2}\right]=\label{eq:k to k+1 z^2-1}
\end{equation}
\[
=\mathbb{E}_{\left\{ k+1|\bm{z}_{0},\bm{z}_{1},\ldots,\bm{z}_{k}\right\} }\left[\|\bm{z}_{k}\|_{2}^{2}-\sum_{l=1}^{P}|\langle\bm{z}_{k},\hat{\bm{u}}_{l}\rangle|^{2}\right]
\]
{\scriptsize{}
\[
\mathbb{=E}_{\left\{ k+1|\bm{z}_{0},\bm{z}_{1},\ldots,\bm{z}_{k}\right\} }\left[\|\bm{z}_{k}\|_{2}^{2}\right]-\mathbb{E}_{\left\{ k+1|\bm{z}_{0},\bm{z}_{1},\ldots,\bm{z}_{k}\right\} }\left[\sum_{l=1}^{P}|\langle\bm{z}_{k},\hat{\bm{u}}_{l}\rangle|^{2}\right].
\]
}We identify the term on the right as:{\footnotesize{}
\[
\mathbb{E}_{\left\{ k+1|\bm{z}_{0},\bm{z}_{1},\ldots,\bm{z}_{k}\right\} }\left[\sum_{l=1}^{P}|\langle\bm{z}_{k},\hat{\bm{u}}_{l}\rangle|^{2}\right]=\mathbb{E}_{\left\{ k+1|\bm{z}_{0},\bm{z}_{1},\ldots,\bm{z}_{k}\right\} }\left[\|U_{P}\bm{\bm{z}}_{k}\|_{2}^{2}\right]
\]
}
\begin{equation}
=\frac{P}{N}\mathbb{\times E}_{\left\{ k+1|\bm{z}_{0},\bm{z}_{1},\ldots,\bm{z}_{k}\right\} }\left[\|\bm{z}_{k}\|_{2}^{2}\right].\label{eq:P/N-1}
\end{equation}
The results from the two equations (\eqref{eq:P/N-1} and \eqref{eq:k to k+1 z^2-1})
above may then be combined to obtain,
\[
\mathbb{E}_{\left\{ k+1|\bm{z}_{0},\bm{z}_{1},\ldots,\bm{z}_{k}\right\} }\left[\|\bm{z}_{k+1}\|_{2}^{2}\right]=
\]
\[
\left(1-\frac{P}{N}\right)\times\mathbb{E}_{\left\{ k|\bm{z}_{0},\bm{z}_{1},\ldots,\bm{z}_{k}-1\right\} }\left[\|\bm{z}_{k}\|_{2}^{2}\right],
\]
 where the expectation on the right hand side includes $k+1\rightarrow k$
accounting for the previous iteration.

Next, apply induction to arrive at the expectation for the whole iterative
sequence up to the $\beta$-th iteration given that $\bm{\bm{z}_{0}}\equiv\bm{x}-\bm{x}_{0}$:
\begin{equation}
\mathbb{E}_{\left\{ \beta+1|\bm{z}_{0}\right\} }\left[\|\bm{z}_{\beta+1}\|_{2}^{2}\right]=\left(1-\frac{P}{N}\right)^{\beta}\|\bm{z}_{0}\|_{2}^{2}\,\,\forall\;\beta\in1,2,3,\ldots.\label{eq:rkos-convergence-term-1}
\end{equation}

\paragraph{Asymptotic Convergence}

The statistical ensemble average of the above Equation \eqref{eq:zl2-1}
for the $\beta$-th iteration yields the convergence result given
in Equation \eqref{eq:rkos-convergence-term-1}. These results assume
random variables identically and independently distributed, but compare
well to others in the literature, such as the convergence result in
Strohmer \cite{Strohmer09}.

The theoretical convergence iterative limit for uniform random IID
sampling was compared to numerical simulations using random solution
vector point on a unit sphere. Equation \eqref{eq:rkos-convergence-2-1}
has an asymptotic form:
\begin{eqnarray}
\frac{\mathbb{E}_{\{\beta+1|\bm{z}_{0}}\}}{\left[\Vert\vct{z}_{\beta+1}\Vert_{2}^{2}\right]}{\Vert\vct{z}_{0}\Vert_{2}^{2}} & =\label{eq:rkos-convergence-2-1}
\end{eqnarray}
\[
\;\lim_{\beta\rightarrow\infty}\left[1-\frac{P}{N}\right]^{\beta}\simeq e^{-\beta P/N}
\]
\[
\; P=\dim(S_{P}),\;\beta\gg1,2,3,\ldots\rightarrow k\in P,2P,3P,\ldots.
\]
For comparison, recall the convergence for RK method of Strohmer
for IID measurements with  $R=N$ is approximately:
\begin{equation}
\frac{\mathbb{E}_{\{k+1|\bm{z}_{0}}\}}{\left[\Vert\vct{z}_{k+1}\Vert_{2}^{2}\right]}{\Vert\vct{z}_{0}\Vert_{2}^{2}}=\left[1-\frac{1}{N}\right]^{k}\label{eq:rkconvergence-2-1}
\end{equation}
\[
\lim_{k\rightarrow\infty}\left[1-\frac{1}{N}\right]^{k}\simeq e^{-k/N}\,\,\forall\; k\gg1,2,3,\ldots.
\]
Estimated noise bound convergence complexity to $\epsilon$ error
is $\bigO(N^{2})$. Since the value of $\bm{z}_{0}$ is given, the
expectation is known to be the same.

\paragraph{Theory and Simulation}

Simulations in reference \cite{Wallace2013} compare theory to Gaussian
IID with noise variance added to the measurements with magnitude
$\beta=0.05$ (about five percent) and iteration termination at $\beta=0.05/4=0.0125$.
In the first problem, the exact solution $\bm{x}$ is chosen as a
random point on the unit sphere - which is illustrated in Figure \ref{Sta:Unit-sphere-2-D-1}.
In a second problem, a measurement of the standard phantom using parallel
beam measurements is included, which contains coherent measurements. 

\subsubsection{QR Representation}

An alternative method for finding the expected convergence of the
RKOS iterative block Kaczmarz method used to solve $A\bm{x}^{*}=\bm{b^{*}}$
for $\left(\bm{x}^{*};\bm{b}^{*}\right)\in\R^{N},$ and $A\in\mathbb{R}^{N\times N}$
is considered below. The formalism is slightly more rigorous and contemporary,
allows direct computation of matrix quantities (instead of recursive
GS), but is consistent with the former method of finding the orthogonal
projection subspaces $U_{i}$. 

The method includes sufficient algebra to allow representation of
the Kaczmarz orthogonal block iterative process subject to the Smith
Solmon Wagner \cite{smith1977} inequality, by incorporating the
subspace projection concepts from Galantai \emph{et a l}\cite{galantai2003projectors}. 

In this work, it is assumed that measurement matrix $A\in\mathbb{R}^{N\times N}$
is square full row rank, however, the results may be extended to cases
where $M\geq N$ with proper modification.

\paragraph{Approach}

The $i$-th block iteration of the RKOS selects blocks of $M_{i}$-rows
of matrix $A$ to form $A_{i}$. In general, the blocks may be selected
to allow overlapping rows or unique row selections per cycle, in natural
row order or via random \emph{a priori }partitioning into the set
$\left\{ M_{1},M_{2},\ldots,M_{k}\right\} $ of row blocks comprising
$A_{i}\in\mathbb{R}^{M_{i}\times N}$. However, in the following analysis,
we assume set is subject to 
\begin{equation}
\sum_{i=1}^{k}M_{i}=N\label{eq:constraint-1}
\end{equation}
which applies to the case in which rows are selected uniquely without
replacement for each cycle.

Let $\H$ be a Hilbert space having a defined inner product and finite
norm. Let the measurement matrix $A\in\mathbb{R}^{N\times N}$ be
full row rank in $\H$ and segmented into $k$-blocks according to
\[
I_{N}=\left[E_{1},\ldots,E_{k}\right]\,\,\left(E_{i}\in\mathbb{R}^{N\times M_{i}},\: i=1,\ldots,k\right)
\]
 where $E_{i}$ is a set of $M_{i}$-column index vectors (which may
be non-contiguous) of the identity matrix $I_{NxN}$ to form $A_{i}^{T}=A^{T}E_{i}.$%
\footnote{To understand the sampling vector $E_{i}$, consider the following
example.

Let $I_{6,6}$ be the identity matrix and select non-continguous sampling
set $M_{i}=\{3,5,6\}$ and form $E_{i}$ as $E_{i}(3,5,6)=\begin{bmatrix}0 & 0 & 0\\
0 & 0 & 0\\
1 & 0 & 0\\
0 & 0 & 0\\
0 & 1 & 0\\
0 & 0 & 1
\end{bmatrix}$ and $E_{i}(3,5,6)^{T}=\begin{bmatrix}0 & 0 & 1 & 0 & 0 & 0\\
0 & 0 & 0 & 0 & 1 & 0\\
0 & 0 & 0 & 0 & 0 & 1
\end{bmatrix}$. %
} 

The segmentation of the blocks and the order of blocks is stationary
with respect to iteration number in this treatment.

\paragraph{QR and Gram-Schmidt}

In the RKOS algorithm, the process of decomposing $A_{i}^{T}$ into
the QR \cite{golub-vanloan:1996,Meyer:2000:MAA:343374} factorization
performs the Gram-Schmidt process for orthogonalization. Algorithm
\eqref{alg:Subspace-Kaczmarz-Projections-1-1} recursively solves
for the orthonormal set and allows recursive computation of the projections
of exact solution $\bm{x}^{*}$ onto the the orthogonal basis in terms
measurements $\bm{b}_{i}.$

Direct QR decomposition for row block $A_{i}$ is noted to be 
\begin{equation}
A_{i}^{T}=A^{T}E_{i}=Q_{i}R_{i}=U_{i}R_{i}\label{eq:QRUR-1}
\end{equation}
 is equivalent to GS and may be directly computed%
\footnote{The transpose is needed since the columns of $U_{i}$ are the rows
of $A_{i}$ block%
}, where $U_{i}\in\mathbb{R}^{N\times M_{i}}$ is the $i$-th orthonormal
basis (columns) constructed from the $M_{i}$-rows randomly selected
from matrix $A$, and $R_{i}\in\mathbb{R}^{M_{i}\times M_{i}}$ is
upper triangular matrix. It is important to note that matrix $Q_{i}=U_{i}$
in the RKOS algorithm \eqref{alg:Subspace-Kaczmarz-Projections-1-1}.

For $\bm{x}\in\H$, define the $j$-th iterative error estimate as,
$\bm{z}_{j}\equiv\bm{\bm{x}^{*}}_{i}-\bm{x}_{j}$ and $\bm{\widetilde{z}}_{j}\equiv\bm{x}_{i}^{*}+\bm{\epsilon}_{x}(i)-\bm{x}_{j}$
respectively without and with noise, where $\bm{x}_{j}$ is the $j$-th
iterative estimate for the $i$-th block projection of $k$-blocks
per cycle; $\bm{x}^{*}$ is the desired noise-free solution to $A\bm{x}^{*}=\bm{b}^{*}$;
$\bm{x}_{i}^{*}$ is the $i$-th block estimate of the noise free
solution; and $\bm{\epsilon}_{x}(i)$ is the $i$-th propagated measurement
noise vector in the current basis.%
\footnote{It should be noted that the noise terms are generally not separable
in practice, but are explicitly shown here in order to facilitate
the analysis.%
}

The simple block Kaczmarz's equation (without noise) using the orthogonal
projection matrix $U_{i}$ may be written as 
\begin{eqnarray}
\bm{x}_{j+1} & = & \bm{x}_{j}+U_{i}U_{i}^{T}(\bm{\bm{x}^{*}}-\bm{x}_{j})\label{eq:RK1-1}\\
\bm{\bm{x}^{*}}-\bm{x}_{j+1} & = & \bm{x}^{*}-\bm{x}_{j}-U_{i}U_{i}^{T}(\bm{\bm{x}^{*}}-\bm{x}_{j})\nonumber \\
\bm{z}_{j+1} & = & \bm{z}_{j}-U_{i}U_{i}^{T}\bm{z}_{j}=(I-U_{i}U_{i}^{T})\bm{z}_{j}\nonumber \\
 &  & \,\,\,(i\equiv j(\mod k)+1).\nonumber 
\end{eqnarray}

In above, notice that $U_{i}$ is orthonormal column matrix, i.e.
$U_{i}^{T}U_{i}=I$ under contraction, but on projection, $U_{i}U_{i}^{T}=P_{i}$
acts to preserves components within the subspace $U_{i}$. The following
relations are noted:
\begin{equation}
A_{i}A_{i}^{T}=(U_{i}R)_{i}^{T}(U_{i}R_{i})=R_{i}^{T}U_{i}^{T}U_{i}R_{i}=R_{i}^{T}R_{i}
\end{equation}
\[
A_{i}^{T}A_{i}=(U_{i}R_{i})(U_{i}R_{i})^{T}=U_{i}R_{i}R_{i}^{T}U_{i}^{T}.
\]
To find the new basis, use definition in Equation \eqref{eq:QRUR-1}
solve to find 
\begin{equation}
U_{i}=A_{i}^{T}R_{i}^{T}(R_{i}R_{i}^{T})^{-1},
\end{equation}
\begin{equation}
U_{i}^{T}U_{i}R_{i}=R_{i}=U_{i}^{T}A_{i}^{T}.
\end{equation}

\paragraph{Block Equations}

Next, consider that measurement vector $\bm{b}$ is comprised of (a)
$\bm{b}^{*}$ the self-consistent error free measurement vector solution
of $A\bm{x}^{*}=\bm{b^{*}}$, and (b) the measurement noise term,
$\bm{\epsilon}_{b}$. Therefore, $A\bm{x}^{*}=\bm{b^{*}}$, $\bm{x}=\bm{x}^{*}+\bm{\epsilon}_{x}$,
$\bm{b}=\bm{b}^{*}+\bm{\epsilon_{b}}$. Then we may find, 
\begin{equation}
A\bm{x}=\bm{b}=\bm{b}^{*}+\bm{\epsilon}_{b},
\end{equation}
\[
A_{i}\bm{x}=\left(U_{i}R_{i}\right)^{T}\bm{x}=R_{i}^{T}U_{i}^{T}\bm{x}=\bm{b}_{i}=\bm{b}_{i}^{*}+\bm{\epsilon}_{b}(i),
\]
where $\bm{b}_{i}=E_{i}\bm{b}$ to obtain the $i$-th under-determined
block estimate for the solution,
\begin{equation}
\bm{x}_{i}=\left(A_{i}^{T}A_{i}\right)^{-1}A_{i}^{T}\bm{b}_{i}=\left(A_{i}^{T}A_{i}\right)^{-1}A_{i}^{T}(\bm{b}_{i}^{*}+\bm{\epsilon}_{b}(i)).
\end{equation}
The next objective is to find the result in the new basis. First,
substitute from Equation \eqref{eq:QRUR-1} and multiply both sides
by $R_{i}$ as follows:
\begin{equation}
A_{i}\bm{x}=(U_{i}R_{i})^{T}\bm{x}=R_{i}^{T}U_{i}^{T}\bm{x}=\bm{b}_{i},
\end{equation}
\[
R_{i}R_{i}^{T}U_{i}^{T}\bm{x}=R_{i}\bm{b}_{i},
\]
\begin{equation}
U_{i}^{T}\left[\bm{x}\right]=(R_{i}R_{i}^{T})^{-1}R_{i}\bm{b}_{i}=(R_{i}R_{i}^{T})^{-1}R_{i}(\bm{b}_{i}^{*}+\bm{\epsilon}_{b}(i)),\label{eq:UiTx-1}
\end{equation}
which has been converted to terms of $U_{i}$ and $R_{i}$. Using
the orthogonality of $U_{i}$, Equation \eqref{eq:UiTx-1} may be
solved for $\bm{x}$ in terms of $U_{i},R_{i}$ as follows:
\begin{eqnarray}
U_{i}U_{i}^{T}\bm{x} & = & U_{i}I\,(R_{i}R_{i}^{T})^{-1}R_{i}\bm{b}_{i}\\
\left(U_{i}U_{i}^{T}\right)\bm{x} & = & \left(U_{i}U_{i}^{T}\right)U_{i}(R_{i}R_{i}^{T})^{-1}R_{i}\bm{b}_{i}\nonumber \\
\rightarrow\bm{x}_{i} & = & \bm{x}_{i}^{*}+\epsilon_{x}(i)=U_{i}(R_{i}R_{i}^{T})^{-1}R_{i}\bm{b}_{i}\nonumber 
\end{eqnarray}
where $\bm{b}_{i}=\bm{b}_{i}^{*}+\bm{\epsilon}_{b}(i)$ and the
contraction of the orthonormal matrix $I=U_{i}^{T}U_{i}$ is used
on the right hand side, $U_{i}U_{i}^{T}$ is non-singular, $\bm{x}=\bm{x}^{*}+\bm{\epsilon}_{x}$,
and the $i$-th block estimate $\bm{x}_{i}=\bm{x}_{i}^{*}+\epsilon_{x}(i)$.
The result $\bm{x}_{i}=U_{i}(R_{i}R_{i}^{T})^{-1}R_{i}\bm{b}_{i}$
may be verified by 
\[
U_{i}^{T}\bm{x}_{i}=U_{i}^{T}U_{i}(R_{i}R_{i}^{T})^{-1}R_{i}\bm{b}_{i}=(R_{i}R_{i}^{T})^{-1}R_{i}\left(\bm{b}_{i}^{*}+\bm{\epsilon}_{b}(i)\right)
\]
 which is equation \eqref{eq:UiTx-1} as expected.

\paragraph{Block Iteration and Noise}

Making the substitutions for the consistent noise free solution $\bm{x}^{*}$
and the measurement noise $\bm{\epsilon}_{b}$, the $j$-th error
difference vector terms are as follows: 
\begin{eqnarray}
\bm{\widetilde{z}}_{j} & \equiv & \bm{z}_{j}+\bm{\epsilon}_{x}=\bm{x}^{*}+\bm{\epsilon}_{x}-\bm{x}_{j},\\
\rightarrow & \bm{\widetilde{z}}_{j}=\bm{x}^{*} & +U_{i}(R_{i}R_{i}^{T})^{-1}R_{i}\bm{\epsilon}_{b}(i)-\bm{x}_{j}.\nonumber 
\end{eqnarray}
The orthogonal block Kaczmarz Equation \eqref{eq:RK1-1} for $\bm{z}_{j+1}\equiv\bm{x}^{*}-\bm{x}_{j+1}$
may be written as follows:
\begin{eqnarray}
\bm{z}_{j+1}+\bm{\epsilon}_{x}(i)= & \left(\bm{z}_{j}+U_{i}(R_{i}R_{i}^{T})^{-1}R_{i}\bm{\epsilon}_{b}(i)\right)\nonumber \\
 & \,\,\,\,-U_{i}U_{i}^{T}\left(\bm{z}_{j}+U_{i}(R_{i}R_{i}^{T})^{-1}R_{i}\bm{\epsilon}_{b}(i)\right)\\
= & (I-U_{i}U_{i}^{T})\left[\bm{z}_{j}+U_{i}(R_{i}R_{i}^{T})^{-1}R_{i}\bm{\epsilon}_{b}(i)\right]\nonumber 
\end{eqnarray}
 or, 
\begin{eqnarray}
\bm{\widetilde{z}}_{j+1} & = & (I-U_{i}U_{i}^{T})\bm{\widetilde{z}}_{j}\,\,(i\equiv j(mod\, k)+1)
\end{eqnarray}
where 
\begin{eqnarray}
\bm{\widetilde{z}}_{j+1} & = & \bm{z}_{j+1}+\bm{\epsilon}_{x}(i)=\bm{x}^{*}+\bm{\epsilon}_{x}(i)-\bm{x}_{j+1},
\end{eqnarray}
and the estimated noise component in the block-row basis is $\bm{\epsilon}_{x}(i)=U_{i}(R_{i}R_{i}^{T})^{-1}R_{i}\bm{\epsilon}_{b}(i)$.
In actual practice, the projected component in the new orthogonal
subspace basis is computed as $U_{i}U_{i}^{T}\bm{x}=U_{i}I\,(R_{i}R_{i}^{T})^{-1}R_{i}\bm{b}_{i}$
from the right hand side, where the value of the under-determined
solution vector $\bm{x}$ for the block estimate is not explicitly
realized.

\paragraph{Cyclical Projections\label{sub:Cyclical-Projections}}

In the notation of Halperin \cite{halperin1962} and Galantai \cite{Galantai2005},
$A_{i}^{T}=A^{T}E_{i}=U_{i}R_{i},$ and the projection operator, null
subspace, and orthonormal condition may be identified as follows:

$P_{M_{j}}=I-U_{j}U_{j}^{T}$, $M_{j}=\mathcal{R}^{\perp}(U_{j}^{T})$,
$U_{j}\in\mathbb{R}^{N\times M_{i}}$ $U_{j}^{T}U_{j}=I_{M_{i}\times M_{i}}$
where during the first cycle, observe that $j=i$ for $j=1,\ldots,k$. 

It is further noted that the cumulative projection and null space
intersection for the $k$-th iteration block are as follows:
\begin{equation}
\Omega=P_{k},\ldots,P_{2}P_{1}=(I-U_{k}U_{k}^{T}),\ldots,(I-U_{1}U_{1}^{T}),
\end{equation}
\begin{equation}
M=\bigcap_{j=1}^{k}\mathcal{R}^{\perp}(U_{j})=\mathcal{R}^{\perp}(\left[U_{1},\ldots,U_{k}\right])=\mathcal{R}^{\perp}(U),\label{eq:nullspace_intersection-1}
\end{equation}
 respectively, with $P_{M}=P_{\mathcal{R}^{\perp}(U)}=I-P_{\mathcal{R}(U)}$.
The Smith Solmon Wagner \cite{smith1977} referenced in Theorem 4
of Galantai \cite{Galantai2005}, has the form
\begin{equation}
\bigg\|\left[(I-U_{k}U_{k}^{T}),\ldots,(I-U_{1}U_{1}^{T})\right]^{N}\bm{z}_{0}-P_{M}\bm{z}_{0}\bigg\|\label{eq:galantai rkos ineq}
\end{equation}
\[
\,\,\,\,\,\,\,\leq c_{SSW}^{N}\bigg\|\bm{z}_{0}-P_{M}\bm{z}_{0}\bigg\|
\]
where $c_{SSW}=\left(\prod_{j=1}^{k-1}\sin^{2}\theta_{j}\right)^{1/2}$and
angle 
\begin{equation}
\theta_{j}=\alpha\left(M_{j},\bigcap_{i=j+1}^{k}M_{i}\right)=\alpha\left(\mathcal{R}^{\perp}(U_{j}),\mathcal{R}^{\perp}(\left[U_{j+1},\ldots,U_{k}\right])\right).
\end{equation}
The above result provides a bound for convergence using linear block
projections

\paragraph{Gram-Schmidt and QR Summary }

The expected statistical convergence method described using Gram-Schmidt
(GS) shows good agreement to experimental simulations. The results
are consistent with Strohmer for $P=1$. The $P$-dimensional orthogonal
subspace method based upon QR gives similar convergence result, and
the deterministic bounds are consistent with the results of Galantai.
In both of the above cases, i.e. Gram-Schmidt and QR decomposition,
the proofs of convergence were based upon IID probability distribution
of the measurement noise and the measurement sampling vectors.

The propagation of measurement noise is seen to be dependent upon
the iterative convergence and general iterative process. An additional
study may be worthwhile to determine a possible method for noise minimization
and feasibility.

\subsubsection{Convergence for Almost Any Probability Distribution}

Although the former methods for RKOS Gram Schmidt and QR assumed IID
random variables, it is noted that application of Theorem \eqref{Halperin}
to Equation \eqref{eq:galantai rkos ineq} in section \eqref{sub:Cyclical-Projections}
yields convergence regardless of the distribution%
\footnote{Note that the span of the solution space must be completely sampled
with non-zero probability %
} of the sampling and IID variates as follows:
\begin{equation}
\lim_{q\to\infty}\left[(I-U_{k}U_{k}^{T}),\ldots,(I-U_{1}U_{1}^{T})\right]^{q}\bm{z}_{0}=P_{M}\bm{z}_{0}.\label{eq:galantai rkos ineq-1}
\end{equation}

As noted before, the block Kaczmarz is an alternating projection method
with $M_{1}=Sp^{\perp}(U_{1}^{T}),\hdots,M_{k}=Sp^{\perp}(U_{k}^{T})$.
Also, $P_{M_{1}}=P_{Sp^{\perp}(U_{1}^{T})},\hdots,P_{M_{k}=Sp^{\perp}(U_{k}^{T})}$
and $M=Sp^{\perp}(U_{1}^{T})\cap\hdots\cap Sp^{\perp}(U_{k}^{T})=Sp^{\perp}(A^{T})$.
Since $A$ is full column rank, $Sp^{\perp}(A^{T})=\{0\}$ and $P_{M}=\{0\}$.
After $q$ cycles, 
\begin{equation}
\bm{z}_{qk}=\bm{x}_{qk}-\bm{x}^{*}=(P_{M_{k}}P_{M_{k-1}}\hdots P_{M_{1}})^{q}(\bm{x}_{0}-\bm{x}^{*}).
\end{equation}
By Theorem~\ref{angles}, $\lim_{q\to\infty}\bm{x}_{qk}-\bm{x}^{*}=0$
and $\lim_{q\to\infty}\bm{x}_{qk}=\bm{x}^{*}$. Here, it should be
noted that orthogonality of $U_{k}$ is consistent with Galantai.

\section{Regular versus Randomized Kaczmarz}

The randomized Kaczmarz's algorithm developed by Strohmer in~\cite{Strohmer2009}
has the following convergence in expectation: 
\begin{equation}
\E\norm{\bm{x}_{qM}-\bm{x}^{*}}_{2}^{2}\leq(1-\frac{1}{\kappa(A)^{2}})^{qM}\norm{\bm{x}_{0}-\bm{x}^{*}}_{2}^{2}\label{eqn:random}
\end{equation}
where $\kappa(A)=\norm{A}_{F}\norm{A^{\dagger}}_{2}$ is the scaled
condition number of matrix $A$ with $A^{\dagger}$ is the left pseudo-inverse
of $A$. The bound for regular Kacmarz is given in Equation~\eqref{eqn:normal}.
Note that we assume $A\in\R^{M\times M}$. Now, we need to compare
$(1-\frac{1}{\norm{A}_{F}^{2}\norm{A^{\dagger}}_{2}^{2}})$ and $(1-\det(AA^{T}))^{1/M}$
to assess which bound is tighter. Let $\sigma_{1}\geq\sigma_{2}\geq\hdots\geq\sigma_{M}>0$
be ordered singular values of $A$. Then, 
\begin{align}
\norm{A^{\dagger}}_{2}^{2} & =1/\sigma_{N}^{2}\\
\norm{A}_{F}^{2} & =\sum_{i=1}^{M}\sigma_{i}^{2}.
\end{align}
Also, note that 
\begin{equation}
AA^{T}=\left[\begin{array}{cccc}
1 & \cos\theta_{12} & \hdots & \cos\theta_{1M}\\
\cos\theta_{21} & 1 & \hdots & \cos\theta_{2M}\\
\vdots & \vdots & \hdots & \vdots\\
\cos\theta_{M1} & \cos\theta_{M2} & \hdots & 1
\end{array}\right]
\end{equation}
where $\theta_{ij}$denotes the angles between the rows $a_{i}$ and
$a_{j}$ of $A$. Then, 
\begin{equation}
\det(AA^{T})\leq\prod_{i=1}^{M}\sum_{j=1}^{M}\cos^{2}\theta_{ij}.
\end{equation}
Note that 
\begin{equation}
\prod_{i=1}^{M}\sigma_{i}^{2}(A)=\prod_{i=1}^{M}\lambda_{i}(A^{T}A)=\det(A^{T}A)=\det(AA^{T})
\end{equation}
therefore 
\begin{equation}
[1-\det(AA^{T})]^{1/M}=(1-\prod_{i=1}^{M}\sigma_{i}^{2})^{1/M}.
\end{equation}
Now, Equations \ref{eqn:random} and \ref{eqn:normal} become: 
\begin{equation}
\E\norm{\bm{x}_{qM}-\bm{x}^{*}}_{2}^{2}\leq(1-\frac{\sigma_{M}^{2}}{\sum_{i=1}^{M}\sigma_{i}^{2}})^{qM}\norm{\bm{x}_{0}-\bm{x}^{*}}_{2}^{2},
\end{equation}
\begin{equation}
\norm{\bm{x}_{qM}-\bm{x}^{*}}_{2}^{2}\leq[(1-\prod_{i=1}^{M}\sigma_{i}^{2})^{1/M}]^{qM}\norm{\bm{x}_{0}-\bm{x}^{*}}_{2}^{2}.
\end{equation}

\section{Experimental Results}

Here, we compare our angle-based randomization with norm-based randomization
of Strohmer \cite{Strohmer2009} in the context of measurement methods.
In particular, a phantom image was used as the solution in simulation
experiments \cite{Herman2009}. Figure \ref{fig:Example-convergence-result fan}
shows that our approach (angle-based randomization) provides a better
convergence rate over the randomized Kaczmarz (norm-based randomization)
in the case of fan-beam sampling. However, our method is computationally
more complex, and therefore we devised another algorithm (explained
in the next following section) that addresses this issue.

The following experiments compare Kaczmarz (K), randomized Kaczmarz
(RK), and randomized Kaczmarz hyperplane angles (RKHA) via simulations.
The objective is to illustrate the effect of row randomization upon
the convergence and observe the dependence upon the sampling methods.
\begin{figure}
\centering{}Angle Probability Distribution for Random Sampling Tomography\vspace{-33pt}
\makebox[1\columnwidth]{%
\subfloat[$\theta(i,j)$ probability density distribution ]{\centering{} % ctangles-compare-3-1.tex
\begin{tikzpicture}[scale=1.][gnuplot]
\path (0.000,0.000) rectangle (8.000,6.000);
\gpfill{rgb color={1.000,1.000,1.000}} (1.196,0.616)--(7.446,0.616)--(7.446,5.074)--(1.196,5.074)--cycle;
\gpcolor{color=gp lt color border}
\gpsetlinetype{gp lt border}
\gpsetlinewidth{.1500}
\draw[gp path] (1.196,0.616)--(1.196,5.074)--(7.446,5.074)--(7.446,0.616)--cycle;
\gpsetlinewidth{0.150}
\draw[gp path] (1.196,0.616)--(1.447,0.616);
\draw[gp path] (7.447,0.616)--(7.196,0.616);
\gpcolor{rgb color={0.000,0.000,0.000}}
\node[gp node right,font={\fontsize{8pt}{9.6pt}\selectfont}] at (1.012,0.616) {$10^{-4}$};
\gpcolor{color=gp lt color border}
\draw[gp path] (1.196,1.287)--(1.321,1.287);
\draw[gp path] (7.447,1.287)--(7.322,1.287);
\draw[gp path] (1.196,1.680)--(1.321,1.680);
\draw[gp path] (7.447,1.680)--(7.322,1.680);
\draw[gp path] (1.196,1.958)--(1.321,1.958);
\draw[gp path] (7.447,1.958)--(7.322,1.958);
\draw[gp path] (1.196,2.174)--(1.321,2.174);
\draw[gp path] (7.447,2.174)--(7.322,2.174);
\draw[gp path] (1.196,2.351)--(1.321,2.351);
\draw[gp path] (7.447,2.351)--(7.322,2.351);
\draw[gp path] (1.196,2.500)--(1.321,2.500);
\draw[gp path] (7.447,2.500)--(7.322,2.500);
\draw[gp path] (1.196,2.629)--(1.321,2.629);
\draw[gp path] (7.447,2.629)--(7.322,2.629);
\draw[gp path] (1.196,2.743)--(1.321,2.743);
\draw[gp path] (7.447,2.743)--(7.322,2.743);
\draw[gp path] (1.196,2.845)--(1.447,2.845);
\draw[gp path] (7.447,2.845)--(7.196,2.845);
\gpcolor{rgb color={0.000,0.000,0.000}}
\node[gp node right,font={\fontsize{8pt}{9.6pt}\selectfont}] at (1.012,2.845) {$10^{-3}$};
\gpcolor{color=gp lt color border}
\draw[gp path] (1.196,3.517)--(1.321,3.517);
\draw[gp path] (7.447,3.517)--(7.322,3.517);
\draw[gp path] (1.196,3.909)--(1.321,3.909);
\draw[gp path] (7.447,3.909)--(7.322,3.909);
\draw[gp path] (1.196,4.188)--(1.321,4.188);
\draw[gp path] (7.447,4.188)--(7.322,4.188);
\draw[gp path] (1.196,4.404)--(1.321,4.404);
\draw[gp path] (7.447,4.404)--(7.322,4.404);
\draw[gp path] (1.196,4.580)--(1.321,4.580);
\draw[gp path] (7.447,4.580)--(7.322,4.580);
\draw[gp path] (1.196,4.730)--(1.321,4.730);
\draw[gp path] (7.447,4.730)--(7.322,4.730);
\draw[gp path] (1.196,4.859)--(1.321,4.859);
\draw[gp path] (7.447,4.859)--(7.322,4.859);
\draw[gp path] (1.196,4.973)--(1.321,4.973);
\draw[gp path] (7.447,4.973)--(7.322,4.973);
\draw[gp path] (1.196,5.075)--(1.447,5.075);
\draw[gp path] (7.447,5.075)--(7.196,5.075);
\gpcolor{rgb color={0.000,0.000,0.000}}
\node[gp node right,font={\fontsize{8pt}{9.6pt}\selectfont}] at (1.012,5.075) {$10^{-2}$};
\gpcolor{color=gp lt color border}
\draw[gp path] (1.196,0.616)--(1.196,0.867);
\draw[gp path] (1.196,5.075)--(1.196,4.824);
\gpcolor{rgb color={0.000,0.000,0.000}}
\node[gp node center,font={\fontsize{8pt}{9.6pt}\selectfont}] at (1.196,0.308) {0};
\gpcolor{color=gp lt color border}
\draw[gp path] (2.238,0.616)--(2.238,0.867);
\draw[gp path] (2.238,5.075)--(2.238,4.824);
\gpcolor{rgb color={0.000,0.000,0.000}}
\node[gp node center,font={\fontsize{8pt}{9.6pt}\selectfont}] at (2.238,0.308) {20};
\gpcolor{color=gp lt color border}
\draw[gp path] (3.280,0.616)--(3.280,0.867);
\draw[gp path] (3.280,5.075)--(3.280,4.824);
\gpcolor{rgb color={0.000,0.000,0.000}}
\node[gp node center,font={\fontsize{8pt}{9.6pt}\selectfont}] at (3.280,0.308) {40};
\gpcolor{color=gp lt color border}
\draw[gp path] (4.322,0.616)--(4.322,0.867);
\draw[gp path] (4.322,5.075)--(4.322,4.824);
\gpcolor{rgb color={0.000,0.000,0.000}}
\node[gp node center,font={\fontsize{8pt}{9.6pt}\selectfont}] at (4.322,0.308) {60};
\gpcolor{color=gp lt color border}
\draw[gp path] (5.363,0.616)--(5.363,0.867);
\draw[gp path] (5.363,5.075)--(5.363,4.824);
\gpcolor{rgb color={0.000,0.000,0.000}}
\node[gp node center,font={\fontsize{8pt}{9.6pt}\selectfont}] at (5.363,0.308) {80};
\gpcolor{color=gp lt color border}
\draw[gp path] (6.405,0.616)--(6.405,0.867);
\draw[gp path] (6.405,5.075)--(6.405,4.824);
\gpcolor{rgb color={0.000,0.000,0.000}}
\node[gp node center,font={\fontsize{8pt}{9.6pt}\selectfont}] at (6.405,0.308) {100};
\gpcolor{color=gp lt color border}
\draw[gp path] (7.447,0.616)--(7.447,0.867);
\draw[gp path] (7.447,5.075)--(7.447,4.824);
\gpcolor{rgb color={0.000,0.000,0.000}}
\node[gp node center,font={\fontsize{8pt}{9.6pt}\selectfont}] at (7.447,0.308) {120};
\gpcolor{color=gp lt color border}
\draw[gp path] (1.196,5.075)--(1.196,0.616)--(7.447,0.616)--(7.447,5.075)--cycle;
\gpcolor{rgb color={0.49,0.73,.91}}
\gpsetlinewidth{2.500}
\gpsetpointsize{8}
\gppoint{gp mark 12}{(1.196,4.404)}
\gppoint{gp mark 12}{(5.010,2.390)}
\gppoint{gp mark 12}{(5.027,2.390)}
\gppoint{gp mark 12}{(5.037,2.390)}
\gppoint{gp mark 12}{(5.069,2.390)}
\gppoint{gp mark 12}{(5.079,2.390)}
\gppoint{gp mark 12}{(5.093,3.062)}
\gppoint{gp mark 12}{(5.110,2.390)}
\gppoint{gp mark 12}{(5.117,2.390)}
\gppoint{gp mark 12}{(5.134,2.390)}
\gppoint{gp mark 12}{(5.141,2.390)}
\gppoint{gp mark 12}{(5.148,2.390)}
\gppoint{gp mark 12}{(5.176,2.390)}
\gppoint{gp mark 12}{(5.187,2.390)}
\gppoint{gp mark 12}{(5.193,2.390)}
\gppoint{gp mark 12}{(5.200,2.390)}
\gppoint{gp mark 12}{(5.204,2.390)}
\gppoint{gp mark 12}{(5.207,2.390)}
\gppoint{gp mark 12}{(5.214,2.390)}
\gppoint{gp mark 12}{(5.228,3.062)}
\gppoint{gp mark 12}{(5.235,2.390)}
\gppoint{gp mark 12}{(5.239,2.390)}
\gppoint{gp mark 12}{(5.242,2.390)}
\gppoint{gp mark 12}{(5.263,3.062)}
\gppoint{gp mark 12}{(5.266,3.733)}
\gppoint{gp mark 12}{(5.277,2.390)}
\gppoint{gp mark 12}{(5.287,2.390)}
\gppoint{gp mark 12}{(5.294,3.062)}
\gppoint{gp mark 12}{(5.301,2.390)}
\gppoint{gp mark 12}{(5.305,2.390)}
\gppoint{gp mark 12}{(5.308,3.454)}
\gppoint{gp mark 12}{(5.315,2.390)}
\gppoint{gp mark 12}{(5.318,2.390)}
\gppoint{gp mark 12}{(5.322,2.390)}
\gppoint{gp mark 12}{(5.325,2.390)}
\gppoint{gp mark 12}{(5.329,3.062)}
\gppoint{gp mark 12}{(5.332,3.062)}
\gppoint{gp mark 12}{(5.343,3.454)}
\gppoint{gp mark 12}{(5.346,3.062)}
\gppoint{gp mark 12}{(5.353,2.390)}
\gppoint{gp mark 12}{(5.357,3.062)}
\gppoint{gp mark 12}{(5.360,3.062)}
\gppoint{gp mark 12}{(5.363,3.454)}
\gppoint{gp mark 12}{(5.367,3.062)}
\gppoint{gp mark 12}{(5.370,3.454)}
\gppoint{gp mark 12}{(5.374,3.454)}
\gppoint{gp mark 12}{(5.377,3.062)}
\gppoint{gp mark 12}{(5.381,3.062)}
\gppoint{gp mark 12}{(5.384,2.390)}
\gppoint{gp mark 12}{(5.388,3.062)}
\gppoint{gp mark 12}{(5.391,3.454)}
\gppoint{gp mark 12}{(5.395,3.454)}
\gppoint{gp mark 12}{(5.398,2.390)}
\gppoint{gp mark 12}{(5.405,3.454)}
\gppoint{gp mark 12}{(5.409,3.733)}
\gppoint{gp mark 12}{(5.412,2.390)}
\gppoint{gp mark 12}{(5.416,3.733)}
\gppoint{gp mark 12}{(5.419,3.733)}
\gppoint{gp mark 12}{(5.422,3.062)}
\gppoint{gp mark 12}{(5.426,3.454)}
\gppoint{gp mark 12}{(5.433,3.062)}
\gppoint{gp mark 12}{(5.443,3.733)}
\gppoint{gp mark 12}{(5.447,4.125)}
\gppoint{gp mark 12}{(5.450,3.949)}
\gppoint{gp mark 12}{(5.454,3.454)}
\gppoint{gp mark 12}{(5.457,3.454)}
\gppoint{gp mark 12}{(5.461,2.390)}
\gppoint{gp mark 12}{(5.464,3.454)}
\gppoint{gp mark 12}{(5.468,2.390)}
\gppoint{gp mark 12}{(5.471,4.404)}
\gppoint{gp mark 12}{(5.475,3.062)}
\gppoint{gp mark 12}{(5.478,4.275)}
\gppoint{gp mark 12}{(5.481,3.454)}
\gppoint{gp mark 12}{(5.485,3.062)}
\gppoint{gp mark 12}{(5.488,2.390)}
\gppoint{gp mark 12}{(5.492,4.404)}
\gppoint{gp mark 12}{(5.495,3.733)}
\gppoint{gp mark 12}{(5.499,3.454)}
\gppoint{gp mark 12}{(5.502,3.062)}
\gppoint{gp mark 12}{(5.506,2.390)}
\gppoint{gp mark 12}{(5.509,3.949)}
\gppoint{gp mark 12}{(5.513,3.733)}
\gppoint{gp mark 12}{(5.516,3.062)}
\gppoint{gp mark 12}{(5.520,3.062)}
\gppoint{gp mark 12}{(5.523,2.390)}
\gppoint{gp mark 12}{(5.527,4.275)}
\gppoint{gp mark 12}{(5.530,3.949)}
\gppoint{gp mark 12}{(5.534,2.390)}
\gppoint{gp mark 12}{(5.537,4.275)}
\gppoint{gp mark 12}{(5.540,3.949)}
\gppoint{gp mark 12}{(5.544,4.125)}
\gppoint{gp mark 12}{(5.547,2.390)}
\gppoint{gp mark 12}{(5.551,4.620)}
\gppoint{gp mark 12}{(5.554,3.949)}
\gppoint{gp mark 12}{(5.558,4.125)}
\gppoint{gp mark 12}{(5.561,4.518)}
\gppoint{gp mark 12}{(5.565,3.949)}
\gppoint{gp mark 12}{(5.568,3.454)}
\gppoint{gp mark 12}{(5.572,4.275)}
\gppoint{gp mark 12}{(5.575,3.733)}
\gppoint{gp mark 12}{(5.579,4.518)}
\gppoint{gp mark 12}{(5.582,3.949)}
\gppoint{gp mark 12}{(5.586,3.949)}
\gppoint{gp mark 12}{(5.589,4.125)}
\gppoint{gp mark 12}{(5.593,4.518)}
\gppoint{gp mark 12}{(5.596,3.733)}
\gppoint{gp mark 12}{(5.599,3.062)}
\gppoint{gp mark 12}{(5.606,4.275)}
\gppoint{gp mark 12}{(5.610,4.518)}
\gppoint{gp mark 12}{(5.613,4.275)}
\gppoint{gp mark 12}{(5.617,3.062)}
\gppoint{gp mark 12}{(5.620,3.733)}
\gppoint{gp mark 12}{(5.624,3.733)}
\gppoint{gp mark 12}{(5.627,4.125)}
\gppoint{gp mark 12}{(5.631,4.125)}
\gppoint{gp mark 12}{(5.634,3.949)}
\gppoint{gp mark 12}{(5.638,4.404)}
\gppoint{gp mark 12}{(5.641,3.733)}
\gppoint{gp mark 12}{(5.645,4.125)}
\gppoint{gp mark 12}{(5.648,3.733)}
\gppoint{gp mark 12}{(5.652,3.733)}
\gppoint{gp mark 12}{(5.655,4.275)}
\gppoint{gp mark 12}{(5.658,3.062)}
\gppoint{gp mark 12}{(5.662,3.454)}
\gppoint{gp mark 12}{(5.665,4.404)}
\gppoint{gp mark 12}{(5.669,4.404)}
\gppoint{gp mark 12}{(5.672,3.733)}
\gppoint{gp mark 12}{(5.676,3.949)}
\gppoint{gp mark 12}{(5.679,3.949)}
\gppoint{gp mark 12}{(5.683,4.125)}
\gppoint{gp mark 12}{(5.686,4.620)}
\gppoint{gp mark 12}{(5.690,4.125)}
\gppoint{gp mark 12}{(5.693,4.125)}
\gppoint{gp mark 12}{(5.697,3.949)}
\gppoint{gp mark 12}{(5.700,3.062)}
\gppoint{gp mark 12}{(5.704,2.390)}
\gppoint{gp mark 12}{(5.707,4.125)}
\gppoint{gp mark 12}{(5.710,4.518)}
\gppoint{gp mark 12}{(5.714,4.404)}
\gppoint{gp mark 12}{(5.717,3.733)}
\gppoint{gp mark 12}{(5.721,3.733)}
\gppoint{gp mark 12}{(5.724,3.949)}
\gppoint{gp mark 12}{(5.728,4.404)}
\gppoint{gp mark 12}{(5.731,3.733)}
\gppoint{gp mark 12}{(5.735,3.733)}
\gppoint{gp mark 12}{(5.738,4.125)}
\gppoint{gp mark 12}{(5.742,3.733)}
\gppoint{gp mark 12}{(5.745,4.125)}
\gppoint{gp mark 12}{(5.749,3.733)}
\gppoint{gp mark 12}{(5.752,4.518)}
\gppoint{gp mark 12}{(5.756,4.125)}
\gppoint{gp mark 12}{(5.759,3.949)}
\gppoint{gp mark 12}{(5.763,3.454)}
\gppoint{gp mark 12}{(5.766,4.404)}
\gppoint{gp mark 12}{(5.769,4.125)}
\gppoint{gp mark 12}{(5.773,4.518)}
\gppoint{gp mark 12}{(5.776,4.404)}
\gppoint{gp mark 12}{(5.780,4.275)}
\gppoint{gp mark 12}{(5.783,4.620)}
\gppoint{gp mark 12}{(5.787,4.712)}
\gppoint{gp mark 12}{(5.790,4.125)}
\gppoint{gp mark 12}{(5.794,4.404)}
\gppoint{gp mark 12}{(5.797,4.518)}
\gppoint{gp mark 12}{(5.801,4.125)}
\gppoint{gp mark 12}{(5.804,4.125)}
\gppoint{gp mark 12}{(5.808,4.275)}
\gppoint{gp mark 12}{(5.811,3.733)}
\gppoint{gp mark 12}{(5.815,4.404)}
\gppoint{gp mark 12}{(5.818,3.733)}
\gppoint{gp mark 12}{(5.822,4.275)}
\gppoint{gp mark 12}{(5.825,4.404)}
\gppoint{gp mark 12}{(5.828,3.949)}
\gppoint{gp mark 12}{(5.832,4.125)}
\gppoint{gp mark 12}{(5.835,4.518)}
\gppoint{gp mark 12}{(5.839,3.733)}
\gppoint{gp mark 12}{(5.842,4.404)}
\gppoint{gp mark 12}{(5.846,4.125)}
\gppoint{gp mark 12}{(5.849,4.404)}
\gppoint{gp mark 12}{(5.853,4.874)}
\gppoint{gp mark 12}{(5.856,3.733)}
\gppoint{gp mark 12}{(5.860,4.404)}
\gppoint{gp mark 12}{(5.863,4.125)}
\gppoint{gp mark 12}{(5.867,4.404)}
\gppoint{gp mark 12}{(5.870,4.518)}
\gppoint{gp mark 12}{(5.874,4.518)}
\gppoint{gp mark 12}{(5.877,4.620)}
\gppoint{gp mark 12}{(5.881,4.404)}
\gppoint{gp mark 12}{(5.884,4.275)}
\gppoint{gp mark 12}{(5.887,4.275)}
\gppoint{gp mark 12}{(5.891,4.874)}
\gppoint{gp mark 12}{(5.894,3.949)}
\gppoint{gp mark 12}{(5.898,5.013)}
\gppoint{gp mark 12}{(5.901,4.125)}
\gppoint{gp mark 12}{(5.905,4.946)}
\gppoint{gp mark 12}{(5.908,4.518)}
\gppoint{gp mark 12}{(5.912,4.275)}
\gppoint{gp mark 12}{(5.915,4.275)}
\gppoint{gp mark 12}{(5.919,4.518)}
\gppoint{gp mark 12}{(5.922,4.874)}
\gppoint{gp mark 12}{(5.926,3.949)}
\gppoint{gp mark 12}{(5.929,3.454)}
\gppoint{gp mark 12}{(5.933,4.518)}
\gppoint{gp mark 12}{(5.936,3.949)}
\gppoint{gp mark 12}{(5.940,4.125)}
\gppoint{gp mark 12}{(5.943,4.712)}
\gppoint{gp mark 12}{(5.946,4.404)}
\gppoint{gp mark 12}{(5.950,4.796)}
\gppoint{gp mark 12}{(5.953,4.518)}
\gppoint{gp mark 12}{(5.957,4.125)}
\gppoint{gp mark 12}{(5.960,3.949)}
\gppoint{gp mark 12}{(5.964,4.275)}
\gppoint{gp mark 12}{(5.967,4.125)}
\gppoint{gp mark 12}{(5.971,3.733)}
\gppoint{gp mark 12}{(5.974,3.949)}
\gppoint{gp mark 12}{(5.978,4.620)}
\gppoint{gp mark 12}{(5.981,4.275)}
\gppoint{gp mark 12}{(5.985,3.062)}
\gppoint{gp mark 12}{(5.988,3.733)}
\gppoint{gp mark 12}{(5.992,4.518)}
\gppoint{gp mark 12}{(5.995,4.404)}
\gppoint{gp mark 12}{(5.999,3.062)}
\gppoint{gp mark 12}{(6.002,4.404)}
\gppoint{gp mark 12}{(6.005,4.125)}
\gppoint{gp mark 12}{(6.009,4.518)}
\gppoint{gp mark 12}{(6.012,4.275)}
\gppoint{gp mark 12}{(6.016,3.062)}
\gppoint{gp mark 12}{(6.019,4.620)}
\gppoint{gp mark 12}{(6.023,4.275)}
\gppoint{gp mark 12}{(6.026,3.733)}
\gppoint{gp mark 12}{(6.030,3.949)}
\gppoint{gp mark 12}{(6.033,4.712)}
\gppoint{gp mark 12}{(6.037,4.620)}
\gppoint{gp mark 12}{(6.040,4.125)}
\gppoint{gp mark 12}{(6.044,3.949)}
\gppoint{gp mark 12}{(6.047,4.125)}
\gppoint{gp mark 12}{(6.051,4.275)}
\gppoint{gp mark 12}{(6.054,3.062)}
\gppoint{gp mark 12}{(6.057,3.733)}
\gppoint{gp mark 12}{(6.061,3.949)}
\gppoint{gp mark 12}{(6.064,4.125)}
\gppoint{gp mark 12}{(6.068,3.733)}
\gppoint{gp mark 12}{(6.071,4.518)}
\gppoint{gp mark 12}{(6.075,4.404)}
\gppoint{gp mark 12}{(6.078,3.454)}
\gppoint{gp mark 12}{(6.082,3.454)}
\gppoint{gp mark 12}{(6.085,4.125)}
\gppoint{gp mark 12}{(6.089,4.275)}
\gppoint{gp mark 12}{(6.092,3.949)}
\gppoint{gp mark 12}{(6.096,4.125)}
\gppoint{gp mark 12}{(6.099,3.733)}
\gppoint{gp mark 12}{(6.103,4.404)}
\gppoint{gp mark 12}{(6.106,3.949)}
\gppoint{gp mark 12}{(6.110,4.275)}
\gppoint{gp mark 12}{(6.113,2.390)}
\gppoint{gp mark 12}{(6.116,4.404)}
\gppoint{gp mark 12}{(6.120,3.733)}
\gppoint{gp mark 12}{(6.123,4.125)}
\gppoint{gp mark 12}{(6.127,3.062)}
\gppoint{gp mark 12}{(6.130,3.949)}
\gppoint{gp mark 12}{(6.134,4.125)}
\gppoint{gp mark 12}{(6.137,4.275)}
\gppoint{gp mark 12}{(6.141,4.275)}
\gppoint{gp mark 12}{(6.144,3.454)}
\gppoint{gp mark 12}{(6.148,4.125)}
\gppoint{gp mark 12}{(6.151,4.125)}
\gppoint{gp mark 12}{(6.155,3.733)}
\gppoint{gp mark 12}{(6.158,3.733)}
\gppoint{gp mark 12}{(6.162,3.949)}
\gppoint{gp mark 12}{(6.165,3.062)}
\gppoint{gp mark 12}{(6.169,3.454)}
\gppoint{gp mark 12}{(6.172,4.125)}
\gppoint{gp mark 12}{(6.175,4.518)}
\gppoint{gp mark 12}{(6.179,3.454)}
\gppoint{gp mark 12}{(6.182,4.125)}
\gppoint{gp mark 12}{(6.186,3.062)}
\gppoint{gp mark 12}{(6.189,3.062)}
\gppoint{gp mark 12}{(6.193,3.454)}
\gppoint{gp mark 12}{(6.196,4.404)}
\gppoint{gp mark 12}{(6.200,3.454)}
\gppoint{gp mark 12}{(6.203,3.733)}
\gppoint{gp mark 12}{(6.207,2.390)}
\gppoint{gp mark 12}{(6.210,3.062)}
\gppoint{gp mark 12}{(6.214,3.454)}
\gppoint{gp mark 12}{(6.217,3.454)}
\gppoint{gp mark 12}{(6.221,2.390)}
\gppoint{gp mark 12}{(6.224,3.454)}
\gppoint{gp mark 12}{(6.228,3.733)}
\gppoint{gp mark 12}{(6.231,2.390)}
\gppoint{gp mark 12}{(6.238,3.454)}
\gppoint{gp mark 12}{(6.241,3.062)}
\gppoint{gp mark 12}{(6.245,3.062)}
\gppoint{gp mark 12}{(6.248,3.454)}
\gppoint{gp mark 12}{(6.252,3.733)}
\gppoint{gp mark 12}{(6.255,3.733)}
\gppoint{gp mark 12}{(6.259,3.062)}
\gppoint{gp mark 12}{(6.266,3.949)}
\gppoint{gp mark 12}{(6.269,4.125)}
\gppoint{gp mark 12}{(6.273,2.390)}
\gppoint{gp mark 12}{(6.276,3.949)}
\gppoint{gp mark 12}{(6.280,3.062)}
\gppoint{gp mark 12}{(6.283,3.949)}
\gppoint{gp mark 12}{(6.287,3.733)}
\gppoint{gp mark 12}{(6.290,3.733)}
\gppoint{gp mark 12}{(6.293,3.062)}
\gppoint{gp mark 12}{(6.297,3.733)}
\gppoint{gp mark 12}{(6.300,3.454)}
\gppoint{gp mark 12}{(6.304,3.733)}
\gppoint{gp mark 12}{(6.307,3.454)}
\gppoint{gp mark 12}{(6.311,3.062)}
\gppoint{gp mark 12}{(6.314,3.062)}
\gppoint{gp mark 12}{(6.318,3.062)}
\gppoint{gp mark 12}{(6.321,3.454)}
\gppoint{gp mark 12}{(6.325,3.062)}
\gppoint{gp mark 12}{(6.328,3.454)}
\gppoint{gp mark 12}{(6.332,3.062)}
\gppoint{gp mark 12}{(6.335,3.733)}
\gppoint{gp mark 12}{(6.339,2.390)}
\gppoint{gp mark 12}{(6.342,3.062)}
\gppoint{gp mark 12}{(6.346,3.062)}
\gppoint{gp mark 12}{(6.352,3.733)}
\gppoint{gp mark 12}{(6.356,2.390)}
\gppoint{gp mark 12}{(6.359,2.390)}
\gppoint{gp mark 12}{(6.363,3.454)}
\gppoint{gp mark 12}{(6.366,3.062)}
\gppoint{gp mark 12}{(6.373,3.062)}
\gppoint{gp mark 12}{(6.377,2.390)}
\gppoint{gp mark 12}{(6.380,3.062)}
\gppoint{gp mark 12}{(6.387,2.390)}
\gppoint{gp mark 12}{(6.391,3.062)}
\gppoint{gp mark 12}{(6.398,4.275)}
\gppoint{gp mark 12}{(6.415,3.733)}
\gppoint{gp mark 12}{(6.418,3.062)}
\gppoint{gp mark 12}{(6.422,2.390)}
\gppoint{gp mark 12}{(6.425,2.390)}
\gppoint{gp mark 12}{(6.429,2.390)}
\gppoint{gp mark 12}{(6.432,2.390)}
\gppoint{gp mark 12}{(6.436,3.062)}
\gppoint{gp mark 12}{(6.446,3.062)}
\gppoint{gp mark 12}{(6.450,3.062)}
\gppoint{gp mark 12}{(6.457,3.733)}
\gppoint{gp mark 12}{(6.460,2.390)}
\gppoint{gp mark 12}{(6.463,3.062)}
\gppoint{gp mark 12}{(6.467,3.062)}
\gppoint{gp mark 12}{(6.470,2.390)}
\gppoint{gp mark 12}{(6.474,3.062)}
\gppoint{gp mark 12}{(6.488,2.390)}
\gppoint{gp mark 12}{(6.491,3.062)}
\gppoint{gp mark 12}{(6.498,3.062)}
\gppoint{gp mark 12}{(6.502,3.062)}
\gppoint{gp mark 12}{(6.509,2.390)}
\gppoint{gp mark 12}{(6.512,3.062)}
\gppoint{gp mark 12}{(6.519,2.390)}
\gppoint{gp mark 12}{(6.522,2.390)}
\gppoint{gp mark 12}{(6.526,3.062)}
\gppoint{gp mark 12}{(6.533,3.062)}
\gppoint{gp mark 12}{(6.536,2.390)}
\gppoint{gp mark 12}{(6.543,3.454)}
\gppoint{gp mark 12}{(6.550,2.390)}
\gppoint{gp mark 12}{(6.557,3.062)}
\gppoint{gp mark 12}{(6.561,2.390)}
\gppoint{gp mark 12}{(6.568,2.390)}
\gppoint{gp mark 12}{(6.571,2.390)}
\gppoint{gp mark 12}{(6.599,3.062)}
\gppoint{gp mark 12}{(6.630,2.390)}
\gppoint{gp mark 12}{(6.637,2.390)}
\gppoint{gp mark 12}{(6.640,3.062)}
\gppoint{gp mark 12}{(6.644,2.390)}
\gppoint{gp mark 12}{(6.647,3.062)}
\gppoint{gp mark 12}{(6.661,2.390)}
\gppoint{gp mark 12}{(6.672,2.390)}
\gppoint{gp mark 12}{(6.682,2.390)}
\gppoint{gp mark 12}{(6.703,2.390)}
\gppoint{gp mark 12}{(6.710,2.390)}
\gppoint{gp mark 12}{(6.717,2.390)}
\gppoint{gp mark 12}{(6.731,2.390)}
\gppoint{gp mark 12}{(6.745,2.390)}
\gpcolor{rgb color={0.000,0.000,0.000}}
\gpdefrectangularnode{gp plot 1}{\pgfpoint{1.196cm}{0.616cm}}{\pgfpoint{7.447cm}{5.075cm}}
\end{tikzpicture}
}\subfloat[Gramian $G(i,j)$ density distribution]{\centering{} 
\begin{tikzpicture}[scale=1.][gnuplot]
\path (0.000,0.000) rectangle (8.000,6.000);
\gpfill{rgb color={1.000,1.000,1.000}} (1.196,0.616)--(7.446,0.616)--(7.446,5.074)--(1.196,5.074)--cycle;
\gpcolor{color=gp lt color border}
\gpsetlinetype{gp lt border}
\gpsetlinewidth{.500}
\draw[gp path] (1.196,0.616)--(1.196,5.074)--(7.446,5.074)--(7.446,0.616)--cycle;
\gpsetlinewidth{0.50}
\draw[gp path] (1.196,0.616)--(1.447,0.616);
\draw[gp path] (7.447,0.616)--(7.196,0.616);
\gpcolor{rgb color={0.000,0.000,0.000}}
\node[gp node right,font={\fontsize{8pt}{9.6pt}\selectfont}] at (1.012,0.616) {$10^{-5}$};
\gpcolor{color=gp lt color border}
\draw[gp path] (1.196,0.952)--(1.321,0.952);
\draw[gp path] (7.447,0.952)--(7.322,0.952);
\draw[gp path] (1.196,1.148)--(1.321,1.148);
\draw[gp path] (7.447,1.148)--(7.322,1.148);
\draw[gp path] (1.196,1.287)--(1.321,1.287);
\draw[gp path] (7.447,1.287)--(7.322,1.287);
\draw[gp path] (1.196,1.395)--(1.321,1.395);
\draw[gp path] (7.447,1.395)--(7.322,1.395);
\draw[gp path] (1.196,1.483)--(1.321,1.483);
\draw[gp path] (7.447,1.483)--(7.322,1.483);
\draw[gp path] (1.196,1.558)--(1.321,1.558);
\draw[gp path] (7.447,1.558)--(7.322,1.558);
\draw[gp path] (1.196,1.623)--(1.321,1.623);
\draw[gp path] (7.447,1.623)--(7.322,1.623);
\draw[gp path] (1.196,1.680)--(1.321,1.680);
\draw[gp path] (7.447,1.680)--(7.322,1.680);
\draw[gp path] (1.196,1.731)--(1.447,1.731);
\draw[gp path] (7.447,1.731)--(7.196,1.731);
\gpcolor{rgb color={0.000,0.000,0.000}}
\node[gp node right,font={\fontsize{8pt}{9.6pt}\selectfont}] at (1.012,1.731) {$10^{-4}$};
\gpcolor{color=gp lt color border}
\draw[gp path] (1.196,2.066)--(1.321,2.066);
\draw[gp path] (7.447,2.066)--(7.322,2.066);
\draw[gp path] (1.196,2.263)--(1.321,2.263);
\draw[gp path] (7.447,2.263)--(7.322,2.263);
\draw[gp path] (1.196,2.402)--(1.321,2.402);
\draw[gp path] (7.447,2.402)--(7.322,2.402);
\draw[gp path] (1.196,2.510)--(1.321,2.510);
\draw[gp path] (7.447,2.510)--(7.322,2.510);
\draw[gp path] (1.196,2.598)--(1.321,2.598);
\draw[gp path] (7.447,2.598)--(7.322,2.598);
\draw[gp path] (1.196,2.673)--(1.321,2.673);
\draw[gp path] (7.447,2.673)--(7.322,2.673);
\draw[gp path] (1.196,2.737)--(1.321,2.737);
\draw[gp path] (7.447,2.737)--(7.322,2.737);
\draw[gp path] (1.196,2.794)--(1.321,2.794);
\draw[gp path] (7.447,2.794)--(7.322,2.794);
\draw[gp path] (1.196,2.846)--(1.447,2.846);
\draw[gp path] (7.447,2.846)--(7.196,2.846);
\gpcolor{rgb color={0.000,0.000,0.000}}
\node[gp node right,font={\fontsize{8pt}{9.6pt}\selectfont}] at (1.012,2.846) {$10^{-3}$};
\gpcolor{color=gp lt color border}
\draw[gp path] (1.196,3.181)--(1.321,3.181);
\draw[gp path] (7.447,3.181)--(7.322,3.181);
\draw[gp path] (1.196,3.377)--(1.321,3.377);
\draw[gp path] (7.447,3.377)--(7.322,3.377);
\draw[gp path] (1.196,3.517)--(1.321,3.517);
\draw[gp path] (7.447,3.517)--(7.322,3.517);
\draw[gp path] (1.196,3.625)--(1.321,3.625);
\draw[gp path] (7.447,3.625)--(7.322,3.625);
\draw[gp path] (1.196,3.713)--(1.321,3.713);
\draw[gp path] (7.447,3.713)--(7.322,3.713);
\draw[gp path] (1.196,3.788)--(1.321,3.788);
\draw[gp path] (7.447,3.788)--(7.322,3.788);
\draw[gp path] (1.196,3.852)--(1.321,3.852);
\draw[gp path] (7.447,3.852)--(7.322,3.852);
\draw[gp path] (1.196,3.909)--(1.321,3.909);
\draw[gp path] (7.447,3.909)--(7.322,3.909);
\draw[gp path] (1.196,3.960)--(1.447,3.960);
\draw[gp path] (7.447,3.960)--(7.196,3.960);
\gpcolor{rgb color={0.000,0.000,0.000}}
\node[gp node right,font={\fontsize{8pt}{9.6pt}\selectfont}] at (1.012,3.960) {$10^{-2}$};
\gpcolor{color=gp lt color border}
\draw[gp path] (1.196,4.296)--(1.321,4.296);
\draw[gp path] (7.447,4.296)--(7.322,4.296);
\draw[gp path] (1.196,4.492)--(1.321,4.492);
\draw[gp path] (7.447,4.492)--(7.322,4.492);
\draw[gp path] (1.196,4.631)--(1.321,4.631);
\draw[gp path] (7.447,4.631)--(7.322,4.631);
\draw[gp path] (1.196,4.739)--(1.321,4.739);
\draw[gp path] (7.447,4.739)--(7.322,4.739);
\draw[gp path] (1.196,4.828)--(1.321,4.828);
\draw[gp path] (7.447,4.828)--(7.322,4.828);
\draw[gp path] (1.196,4.902)--(1.321,4.902);
\draw[gp path] (7.447,4.902)--(7.322,4.902);
\draw[gp path] (1.196,4.967)--(1.321,4.967);
\draw[gp path] (7.447,4.967)--(7.322,4.967);
\draw[gp path] (1.196,5.024)--(1.321,5.024);
\draw[gp path] (7.447,5.024)--(7.322,5.024);
\draw[gp path] (1.196,5.075)--(1.447,5.075);
\draw[gp path] (7.447,5.075)--(7.196,5.075);
\gpcolor{rgb color={0.000,0.000,0.000}}
\node[gp node right,font={\fontsize{8pt}{9.6pt}\selectfont}] at (1.012,5.075) {$10^{-1}$};
\gpcolor{color=gp lt color border}
\draw[gp path] (1.196,0.616)--(1.196,0.867);
\draw[gp path] (1.196,5.075)--(1.196,4.824);
\gpcolor{rgb color={0.000,0.000,0.000}}
\node[gp node center,font={\fontsize{8pt}{9.6pt}\selectfont}] at (1.196,0.308) {-0.4};
\gpcolor{color=gp lt color border}
\draw[gp path] (1.977,0.616)--(1.977,0.867);
\draw[gp path] (1.977,5.075)--(1.977,4.824);
\gpcolor{rgb color={0.000,0.000,0.000}}
\node[gp node center,font={\fontsize{8pt}{9.6pt}\selectfont}] at (1.977,0.308) {-0.3};
\gpcolor{color=gp lt color border}
\draw[gp path] (2.759,0.616)--(2.759,0.867);
\draw[gp path] (2.759,5.075)--(2.759,4.824);
\gpcolor{rgb color={0.000,0.000,0.000}}
\node[gp node center,font={\fontsize{8pt}{9.6pt}\selectfont}] at (2.759,0.308) {-0.2};
\gpcolor{color=gp lt color border}
\draw[gp path] (3.540,0.616)--(3.540,0.867);
\draw[gp path] (3.540,5.075)--(3.540,4.824);
\gpcolor{rgb color={0.000,0.000,0.000}}
\node[gp node center,font={\fontsize{8pt}{9.6pt}\selectfont}] at (3.540,0.308) {-0.1};
\gpcolor{color=gp lt color border}
\draw[gp path] (4.322,0.616)--(4.322,0.867);
\draw[gp path] (4.322,5.075)--(4.322,4.824);
\gpcolor{rgb color={0.000,0.000,0.000}}
\node[gp node center,font={\fontsize{8pt}{9.6pt}\selectfont}] at (4.322,0.308) {0};
\gpcolor{color=gp lt color border}
\draw[gp path] (5.103,0.616)--(5.103,0.867);
\draw[gp path] (5.103,5.075)--(5.103,4.824);
\gpcolor{rgb color={0.000,0.000,0.000}}
\node[gp node center,font={\fontsize{8pt}{9.6pt}\selectfont}] at (5.103,0.308) {0.1};
\gpcolor{color=gp lt color border}
\draw[gp path] (5.884,0.616)--(5.884,0.867);
\draw[gp path] (5.884,5.075)--(5.884,4.824);
\gpcolor{rgb color={0.000,0.000,0.000}}
\node[gp node center,font={\fontsize{8pt}{9.6pt}\selectfont}] at (5.884,0.308) {0.2};
\gpcolor{color=gp lt color border}
\draw[gp path] (6.666,0.616)--(6.666,0.867);
\draw[gp path] (6.666,5.075)--(6.666,4.824);
\gpcolor{rgb color={0.000,0.000,0.000}}
\node[gp node center,font={\fontsize{8pt}{9.6pt}\selectfont}] at (6.666,0.308) {0.3};
\gpcolor{color=gp lt color border}
\draw[gp path] (7.447,0.616)--(7.447,0.867);
\draw[gp path] (7.447,5.075)--(7.447,4.824);
\gpcolor{rgb color={0.000,0.000,0.000}}
\node[gp node center,font={\fontsize{8pt}{9.6pt}\selectfont}] at (7.447,0.308) {0.4};
\gpcolor{color=gp lt color border}
\draw[gp path] (1.196,5.075)--(1.196,0.616)--(7.447,0.616)--(7.447,5.075)--cycle;
\gpcolor{rgb color={0.70,0.75,.71}}
\gpsetlinewidth{2.500}
\gpsetpointsize{8}
\gppoint{gp mark 12}{(1.507,1.731)}
\gppoint{gp mark 12}{(1.625,1.395)}
\gppoint{gp mark 12}{(1.684,1.927)}
\gppoint{gp mark 12}{(1.743,1.731)}
\gppoint{gp mark 12}{(1.803,1.731)}
\gppoint{gp mark 12}{(1.862,1.927)}
\gppoint{gp mark 12}{(1.921,2.066)}
\gppoint{gp mark 12}{(1.980,2.066)}
\gppoint{gp mark 12}{(2.040,2.174)}
\gppoint{gp mark 12}{(2.099,2.174)}
\gppoint{gp mark 12}{(2.158,2.598)}
\gppoint{gp mark 12}{(2.217,2.637)}
\gppoint{gp mark 12}{(2.276,2.846)}
\gppoint{gp mark 12}{(2.336,2.913)}
\gppoint{gp mark 12}{(2.395,3.073)}
\gppoint{gp mark 12}{(2.454,2.954)}
\gppoint{gp mark 12}{(2.513,3.143)}
\gppoint{gp mark 12}{(2.572,3.299)}
\gppoint{gp mark 12}{(2.632,3.424)}
\gppoint{gp mark 12}{(2.691,3.459)}
\gppoint{gp mark 12}{(2.750,3.511)}
\gppoint{gp mark 12}{(2.809,3.579)}
\gppoint{gp mark 12}{(2.868,3.666)}
\gppoint{gp mark 12}{(2.928,3.774)}
\gppoint{gp mark 12}{(2.987,3.846)}
\gppoint{gp mark 12}{(3.046,3.870)}
\gppoint{gp mark 12}{(3.105,3.909)}
\gppoint{gp mark 12}{(3.164,3.975)}
\gppoint{gp mark 12}{(3.224,4.068)}
\gppoint{gp mark 12}{(3.283,4.093)}
\gppoint{gp mark 12}{(3.342,4.135)}
\gppoint{gp mark 12}{(3.401,4.106)}
\gppoint{gp mark 12}{(3.460,4.191)}
\gppoint{gp mark 12}{(3.520,4.211)}
\gppoint{gp mark 12}{(3.579,4.325)}
\gppoint{gp mark 12}{(3.638,4.314)}
\gppoint{gp mark 12}{(3.697,4.309)}
\gppoint{gp mark 12}{(3.756,4.351)}
\gppoint{gp mark 12}{(3.816,4.393)}
\gppoint{gp mark 12}{(3.875,4.380)}
\gppoint{gp mark 12}{(3.934,4.439)}
\gppoint{gp mark 12}{(3.993,4.429)}
\gppoint{gp mark 12}{(4.052,4.460)}
\gppoint{gp mark 12}{(4.112,4.459)}
\gppoint{gp mark 12}{(4.171,4.466)}
\gppoint{gp mark 12}{(4.230,4.500)}
\gppoint{gp mark 12}{(4.289,4.512)}
\gppoint{gp mark 12}{(4.349,4.561)}
\gppoint{gp mark 12}{(4.408,4.506)}
\gppoint{gp mark 12}{(4.467,4.451)}
\gppoint{gp mark 12}{(4.526,4.499)}
\gppoint{gp mark 12}{(4.585,4.510)}
\gppoint{gp mark 12}{(4.645,4.426)}
\gppoint{gp mark 12}{(4.704,4.390)}
\gppoint{gp mark 12}{(4.763,4.423)}
\gppoint{gp mark 12}{(4.822,4.414)}
\gppoint{gp mark 12}{(4.881,4.350)}
\gppoint{gp mark 12}{(4.941,4.341)}
\gppoint{gp mark 12}{(5.000,4.310)}
\gppoint{gp mark 12}{(5.059,4.279)}
\gppoint{gp mark 12}{(5.118,4.235)}
\gppoint{gp mark 12}{(5.177,4.239)}
\gppoint{gp mark 12}{(5.237,4.197)}
\gppoint{gp mark 12}{(5.296,4.114)}
\gppoint{gp mark 12}{(5.355,4.102)}
\gppoint{gp mark 12}{(5.414,4.057)}
\gppoint{gp mark 12}{(5.473,3.960)}
\gppoint{gp mark 12}{(5.533,3.912)}
\gppoint{gp mark 12}{(5.592,3.814)}
\gppoint{gp mark 12}{(5.651,3.794)}
\gppoint{gp mark 12}{(5.710,3.725)}
\gppoint{gp mark 12}{(5.769,3.675)}
\gppoint{gp mark 12}{(5.829,3.666)}
\gppoint{gp mark 12}{(5.888,3.472)}
\gppoint{gp mark 12}{(5.947,3.308)}
\gppoint{gp mark 12}{(6.006,3.401)}
\gppoint{gp mark 12}{(6.065,3.335)}
\gppoint{gp mark 12}{(6.125,3.227)}
\gppoint{gp mark 12}{(6.184,3.102)}
\gppoint{gp mark 12}{(6.243,3.169)}
\gppoint{gp mark 12}{(6.302,2.821)}
\gppoint{gp mark 12}{(6.361,2.737)}
\gppoint{gp mark 12}{(6.421,2.767)}
\gppoint{gp mark 12}{(6.480,2.598)}
\gppoint{gp mark 12}{(6.539,2.174)}
\gppoint{gp mark 12}{(6.598,2.174)}
\gppoint{gp mark 12}{(6.657,2.066)}
\gppoint{gp mark 12}{(6.717,2.174)}
\gppoint{gp mark 12}{(6.776,1.927)}
\gppoint{gp mark 12}{(6.835,2.337)}
\gppoint{gp mark 12}{(6.954,1.927)}
\gppoint{gp mark 12}{(7.013,1.395)}
\gppoint{gp mark 12}{(7.072,1.395)}
\gppoint{gp mark 12}{(7.131,1.395)}
\gppoint{gp mark 12}{(7.368,1.395)}
\gpcolor{rgb color={0.000,0.000,0.000}}
\gpdefrectangularnode{gp plot 1}{\pgfpoint{1.196cm}{0.616cm}}{\pgfpoint{7.447cm}{5.075cm}}
\end{tikzpicture}
}%
}\protect\caption{\label{fig:random-1-2-1-1}(a) Example angles distribution (y-axis)
from $AA^{T}$ where $\theta_{i,j}=\cos^{-1}(\langle\hat{\bm{a}}_{i},\hat{\bm{a}}_{j}\rangle)\,\,\forall\,\, i,j\,\in\{1,\ldots,M\}$
vs angles (x-axis) degrees using random data acquisition strategy,
(b) Gramian matrix $\langle\hat{\bm{a}}_{i},\hat{\bm{a}}_{j}\rangle$
distribution}
\end{figure}
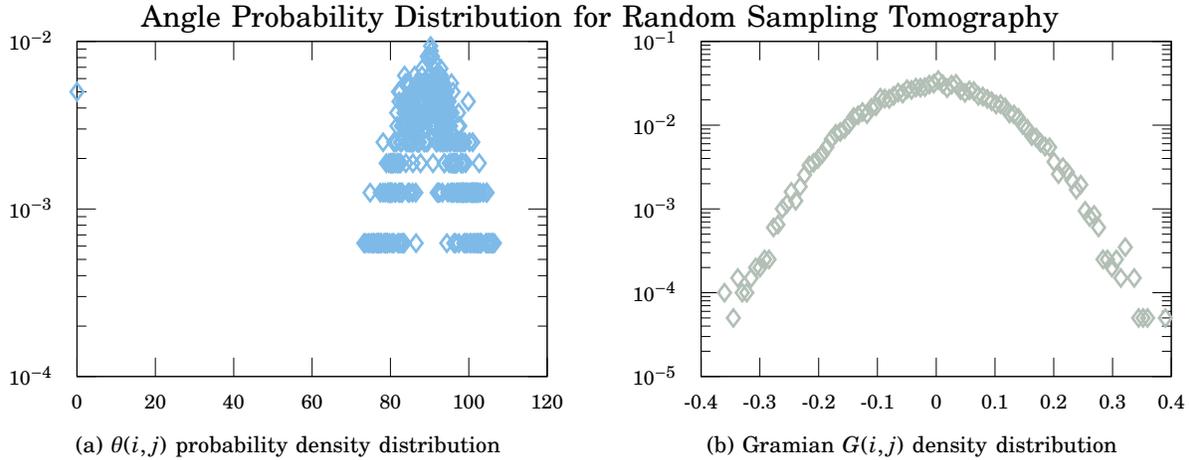
\begin{figure}
\noindent \centering{}Angle Probability Distribution for Fan-Beam
Sampling Tomography\vspace{-33pt}
\makebox[1\columnwidth]{%
\subfloat[$\theta(i,j)$ probability density distribution ]{\centering{} % ctangles-compare-2-1-dblue.tex
\begin{tikzpicture}[scale=1.][gnuplot]
\path (0.000,0.000) rectangle (8.000,6.000);
\gpfill{rgb color={1.000,1.000,1.000}} (1.196,0.616)--(7.446,0.616)--(7.446,5.074)--(1.196,5.074)--cycle;
\gpcolor{color=gp lt color border}
\gpsetlinetype{gp lt border}
\gpsetlinewidth{.500}
\draw[gp path] (1.196,0.616)--(1.196,5.074)--(7.446,5.074)--(7.446,0.616)--cycle;
\gpsetlinewidth{0.50}
\draw[gp path] (1.196,0.616)--(1.447,0.616);
\draw[gp path] (7.447,0.616)--(7.196,0.616);
\gpcolor{rgb color={0.000,0.000,0.000}}
\node[gp node right,font={\fontsize{8pt}{9.6pt}\selectfont}] at (1.012,0.616) {$10^{-4}$};
\gpcolor{color=gp lt color border}
\draw[gp path] (1.196,0.952)--(1.321,0.952);
\draw[gp path] (7.447,0.952)--(7.322,0.952);
\draw[gp path] (1.196,1.148)--(1.321,1.148);
\draw[gp path] (7.447,1.148)--(7.322,1.148);
\draw[gp path] (1.196,1.287)--(1.321,1.287);
\draw[gp path] (7.447,1.287)--(7.322,1.287);
\draw[gp path] (1.196,1.395)--(1.321,1.395);
\draw[gp path] (7.447,1.395)--(7.322,1.395);
\draw[gp path] (1.196,1.483)--(1.321,1.483);
\draw[gp path] (7.447,1.483)--(7.322,1.483);
\draw[gp path] (1.196,1.558)--(1.321,1.558);
\draw[gp path] (7.447,1.558)--(7.322,1.558);
\draw[gp path] (1.196,1.623)--(1.321,1.623);
\draw[gp path] (7.447,1.623)--(7.322,1.623);
\draw[gp path] (1.196,1.680)--(1.321,1.680);
\draw[gp path] (7.447,1.680)--(7.322,1.680);
\draw[gp path] (1.196,1.731)--(1.447,1.731);
\draw[gp path] (7.447,1.731)--(7.196,1.731);
\gpcolor{rgb color={0.000,0.000,0.000}}
\node[gp node right,font={\fontsize{8pt}{9.6pt}\selectfont}] at (1.012,1.731) {$10^{-3}$};
\gpcolor{color=gp lt color border}
\draw[gp path] (1.196,2.066)--(1.321,2.066);
\draw[gp path] (7.447,2.066)--(7.322,2.066);
\draw[gp path] (1.196,2.263)--(1.321,2.263);
\draw[gp path] (7.447,2.263)--(7.322,2.263);
\draw[gp path] (1.196,2.402)--(1.321,2.402);
\draw[gp path] (7.447,2.402)--(7.322,2.402);
\draw[gp path] (1.196,2.510)--(1.321,2.510);
\draw[gp path] (7.447,2.510)--(7.322,2.510);
\draw[gp path] (1.196,2.598)--(1.321,2.598);
\draw[gp path] (7.447,2.598)--(7.322,2.598);
\draw[gp path] (1.196,2.673)--(1.321,2.673);
\draw[gp path] (7.447,2.673)--(7.322,2.673);
\draw[gp path] (1.196,2.737)--(1.321,2.737);
\draw[gp path] (7.447,2.737)--(7.322,2.737);
\draw[gp path] (1.196,2.794)--(1.321,2.794);
\draw[gp path] (7.447,2.794)--(7.322,2.794);
\draw[gp path] (1.196,2.845)--(1.447,2.845);
\draw[gp path] (7.447,2.845)--(7.196,2.845);
\gpcolor{rgb color={0.000,0.000,0.000}}
\node[gp node right,font={\fontsize{8pt}{9.6pt}\selectfont}] at (1.012,2.845) {$10^{-2}$};
\gpcolor{color=gp lt color border}
\draw[gp path] (1.196,3.181)--(1.321,3.181);
\draw[gp path] (7.447,3.181)--(7.322,3.181);
\draw[gp path] (1.196,3.377)--(1.321,3.377);
\draw[gp path] (7.447,3.377)--(7.322,3.377);
\draw[gp path] (1.196,3.517)--(1.321,3.517);
\draw[gp path] (7.447,3.517)--(7.322,3.517);
\draw[gp path] (1.196,3.625)--(1.321,3.625);
\draw[gp path] (7.447,3.625)--(7.322,3.625);
\draw[gp path] (1.196,3.713)--(1.321,3.713);
\draw[gp path] (7.447,3.713)--(7.322,3.713);
\draw[gp path] (1.196,3.788)--(1.321,3.788);
\draw[gp path] (7.447,3.788)--(7.322,3.788);
\draw[gp path] (1.196,3.852)--(1.321,3.852);
\draw[gp path] (7.447,3.852)--(7.322,3.852);
\draw[gp path] (1.196,3.909)--(1.321,3.909);
\draw[gp path] (7.447,3.909)--(7.322,3.909);
\draw[gp path] (1.196,3.960)--(1.447,3.960);
\draw[gp path] (7.447,3.960)--(7.196,3.960);
\gpcolor{rgb color={0.000,0.000,0.000}}
\node[gp node right,font={\fontsize{8pt}{9.6pt}\selectfont}] at (1.012,3.960) {$10^{-1}$};
\gpcolor{color=gp lt color border}
\draw[gp path] (1.196,4.296)--(1.321,4.296);
\draw[gp path] (7.447,4.296)--(7.322,4.296);
\draw[gp path] (1.196,4.492)--(1.321,4.492);
\draw[gp path] (7.447,4.492)--(7.322,4.492);
\draw[gp path] (1.196,4.631)--(1.321,4.631);
\draw[gp path] (7.447,4.631)--(7.322,4.631);
\draw[gp path] (1.196,4.739)--(1.321,4.739);
\draw[gp path] (7.447,4.739)--(7.322,4.739);
\draw[gp path] (1.196,4.828)--(1.321,4.828);
\draw[gp path] (7.447,4.828)--(7.322,4.828);
\draw[gp path] (1.196,4.902)--(1.321,4.902);
\draw[gp path] (7.447,4.902)--(7.322,4.902);
\draw[gp path] (1.196,4.967)--(1.321,4.967);
\draw[gp path] (7.447,4.967)--(7.322,4.967);
\draw[gp path] (1.196,5.024)--(1.321,5.024);
\draw[gp path] (7.447,5.024)--(7.322,5.024);
\draw[gp path] (1.196,5.075)--(1.447,5.075);
\draw[gp path] (7.447,5.075)--(7.196,5.075);
\gpcolor{rgb color={0.000,0.000,0.000}}
\node[gp node right,font={\fontsize{8pt}{9.6pt}\selectfont}] at (1.012,5.075) {$10^{0}$};
\gpcolor{color=gp lt color border}
\draw[gp path] (1.196,0.616)--(1.196,0.867);
\draw[gp path] (1.196,5.075)--(1.196,4.824);
\gpcolor{rgb color={0.000,0.000,0.000}}
\node[gp node center,font={\fontsize{8pt}{9.6pt}\selectfont}] at (1.196,0.308) {0};
\gpcolor{color=gp lt color border}
\draw[gp path] (2.446,0.616)--(2.446,0.867);
\draw[gp path] (2.446,5.075)--(2.446,4.824);
\gpcolor{rgb color={0.000,0.000,0.000}}
\node[gp node center,font={\fontsize{8pt}{9.6pt}\selectfont}] at (2.446,0.308) {20};
\gpcolor{color=gp lt color border}
\draw[gp path] (3.696,0.616)--(3.696,0.867);
\draw[gp path] (3.696,5.075)--(3.696,4.824);
\gpcolor{rgb color={0.000,0.000,0.000}}
\node[gp node center,font={\fontsize{8pt}{9.6pt}\selectfont}] at (3.696,0.308) {40};
\gpcolor{color=gp lt color border}
\draw[gp path] (4.947,0.616)--(4.947,0.867);
\draw[gp path] (4.947,5.075)--(4.947,4.824);
\gpcolor{rgb color={0.000,0.000,0.000}}
\node[gp node center,font={\fontsize{8pt}{9.6pt}\selectfont}] at (4.947,0.308) {60};
\gpcolor{color=gp lt color border}
\draw[gp path] (6.197,0.616)--(6.197,0.867);
\draw[gp path] (6.197,5.075)--(6.197,4.824);
\gpcolor{rgb color={0.000,0.000,0.000}}
\node[gp node center,font={\fontsize{8pt}{9.6pt}\selectfont}] at (6.197,0.308) {80};
\gpcolor{color=gp lt color border}
\draw[gp path] (7.447,0.616)--(7.447,0.867);
\draw[gp path] (7.447,5.075)--(7.447,4.824);
\gpcolor{rgb color={0.000,0.000,0.000}}
\node[gp node center,font={\fontsize{8pt}{9.6pt}\selectfont}] at (7.447,0.308) {100};
\gpcolor{color=gp lt color border}
\draw[gp path] (1.196,5.075)--(1.196,0.616)--(7.447,0.616)--(7.447,5.075)--cycle;
\gpcolor{rgb color={0.49,0.73,.91}}
\gpsetlinewidth{2.500}
\gpsetpointsize{8}
\gppoint{gp mark 12}{(1.196,2.232)}
\gppoint{gp mark 12}{(1.890,1.453)}
\gppoint{gp mark 12}{(1.921,1.453)}
\gppoint{gp mark 12}{(2.035,1.788)}
\gppoint{gp mark 12}{(2.359,1.453)}
\gppoint{gp mark 12}{(2.393,1.788)}
\gppoint{gp mark 12}{(2.406,1.788)}
\gppoint{gp mark 12}{(2.767,1.453)}
\gppoint{gp mark 12}{(3.125,1.453)}
\gppoint{gp mark 12}{(3.227,1.453)}
\gppoint{gp mark 12}{(3.391,1.453)}
\gppoint{gp mark 12}{(3.683,1.453)}
\gppoint{gp mark 12}{(3.787,1.453)}
\gppoint{gp mark 12}{(3.809,1.453)}
\gppoint{gp mark 12}{(4.094,1.453)}
\gppoint{gp mark 12}{(4.243,1.453)}
\gppoint{gp mark 12}{(4.297,1.788)}
\gppoint{gp mark 12}{(4.402,1.453)}
\gppoint{gp mark 12}{(4.535,1.453)}
\gppoint{gp mark 12}{(4.738,1.453)}
\gppoint{gp mark 12}{(4.779,1.453)}
\gppoint{gp mark 12}{(4.788,1.788)}
\gppoint{gp mark 12}{(4.845,1.453)}
\gppoint{gp mark 12}{(4.947,1.453)}
\gppoint{gp mark 12}{(4.969,1.453)}
\gppoint{gp mark 12}{(5.194,1.453)}
\gppoint{gp mark 12}{(5.327,1.453)}
\gppoint{gp mark 12}{(5.365,1.453)}
\gppoint{gp mark 12}{(5.425,1.453)}
\gppoint{gp mark 12}{(5.457,1.453)}
\gppoint{gp mark 12}{(5.460,1.985)}
\gppoint{gp mark 12}{(5.482,1.453)}
\gppoint{gp mark 12}{(5.495,1.453)}
\gppoint{gp mark 12}{(5.498,1.453)}
\gppoint{gp mark 12}{(5.507,1.453)}
\gppoint{gp mark 12}{(5.520,1.453)}
\gppoint{gp mark 12}{(5.539,1.453)}
\gppoint{gp mark 12}{(5.574,1.453)}
\gppoint{gp mark 12}{(5.580,1.453)}
\gppoint{gp mark 12}{(5.609,2.459)}
\gppoint{gp mark 12}{(5.612,1.788)}
\gppoint{gp mark 12}{(5.615,1.985)}
\gppoint{gp mark 12}{(5.618,1.788)}
\gppoint{gp mark 12}{(5.644,2.124)}
\gppoint{gp mark 12}{(5.647,1.453)}
\gppoint{gp mark 12}{(5.653,1.453)}
\gppoint{gp mark 12}{(5.682,1.453)}
\gppoint{gp mark 12}{(5.685,1.453)}
\gppoint{gp mark 12}{(5.773,1.453)}
\gppoint{gp mark 12}{(5.796,1.453)}
\gppoint{gp mark 12}{(5.799,1.788)}
\gppoint{gp mark 12}{(5.808,1.453)}
\gppoint{gp mark 12}{(5.811,1.453)}
\gppoint{gp mark 12}{(5.910,1.985)}
\gppoint{gp mark 12}{(5.922,1.453)}
\gppoint{gp mark 12}{(5.929,1.453)}
\gppoint{gp mark 12}{(5.963,1.453)}
\gppoint{gp mark 12}{(5.970,1.453)}
\gppoint{gp mark 12}{(6.005,1.788)}
\gppoint{gp mark 12}{(6.008,1.453)}
\gppoint{gp mark 12}{(6.030,2.124)}
\gppoint{gp mark 12}{(6.043,1.453)}
\gppoint{gp mark 12}{(6.046,1.453)}
\gppoint{gp mark 12}{(6.049,1.788)}
\gppoint{gp mark 12}{(6.052,1.453)}
\gppoint{gp mark 12}{(6.055,1.453)}
\gppoint{gp mark 12}{(6.058,1.453)}
\gppoint{gp mark 12}{(6.062,1.788)}
\gppoint{gp mark 12}{(6.065,1.788)}
\gppoint{gp mark 12}{(6.077,1.453)}
\gppoint{gp mark 12}{(6.084,1.453)}
\gppoint{gp mark 12}{(6.106,1.453)}
\gppoint{gp mark 12}{(6.125,1.453)}
\gppoint{gp mark 12}{(6.128,1.788)}
\gppoint{gp mark 12}{(6.135,1.788)}
\gppoint{gp mark 12}{(6.141,1.788)}
\gppoint{gp mark 12}{(6.144,1.788)}
\gppoint{gp mark 12}{(6.147,1.453)}
\gppoint{gp mark 12}{(6.150,1.453)}
\gppoint{gp mark 12}{(6.154,1.453)}
\gppoint{gp mark 12}{(6.176,1.453)}
\gppoint{gp mark 12}{(6.179,1.453)}
\gppoint{gp mark 12}{(6.185,1.453)}
\gppoint{gp mark 12}{(6.188,1.453)}
\gppoint{gp mark 12}{(6.195,1.788)}
\gppoint{gp mark 12}{(6.204,1.788)}
\gppoint{gp mark 12}{(6.207,1.453)}
\gppoint{gp mark 12}{(6.211,1.985)}
\gppoint{gp mark 12}{(6.214,1.788)}
\gppoint{gp mark 12}{(6.217,1.788)}
\gppoint{gp mark 12}{(6.223,1.985)}
\gppoint{gp mark 12}{(6.226,2.124)}
\gppoint{gp mark 12}{(6.230,1.453)}
\gppoint{gp mark 12}{(6.233,1.788)}
\gppoint{gp mark 12}{(6.242,2.232)}
\gppoint{gp mark 12}{(6.255,1.453)}
\gppoint{gp mark 12}{(6.258,2.395)}
\gppoint{gp mark 12}{(6.261,1.453)}
\gppoint{gp mark 12}{(6.264,1.453)}
\gppoint{gp mark 12}{(6.268,1.788)}
\gppoint{gp mark 12}{(6.271,1.788)}
\gppoint{gp mark 12}{(6.274,1.788)}
\gppoint{gp mark 12}{(6.277,2.124)}
\gppoint{gp mark 12}{(6.280,1.453)}
\gppoint{gp mark 12}{(6.287,1.788)}
\gppoint{gp mark 12}{(6.290,2.124)}
\gppoint{gp mark 12}{(6.293,2.124)}
\gppoint{gp mark 12}{(6.296,1.453)}
\gppoint{gp mark 12}{(6.299,1.453)}
\gppoint{gp mark 12}{(6.302,1.453)}
\gppoint{gp mark 12}{(6.315,1.453)}
\gppoint{gp mark 12}{(6.318,2.232)}
\gppoint{gp mark 12}{(6.325,1.453)}
\gppoint{gp mark 12}{(6.328,1.453)}
\gppoint{gp mark 12}{(6.331,2.320)}
\gppoint{gp mark 12}{(6.340,1.788)}
\gppoint{gp mark 12}{(6.344,1.788)}
\gppoint{gp mark 12}{(6.347,1.453)}
\gppoint{gp mark 12}{(6.350,1.453)}
\gppoint{gp mark 12}{(6.353,2.459)}
\gppoint{gp mark 12}{(6.356,2.124)}
\gppoint{gp mark 12}{(6.359,1.788)}
\gppoint{gp mark 12}{(6.366,2.232)}
\gppoint{gp mark 12}{(6.369,2.124)}
\gppoint{gp mark 12}{(6.372,1.985)}
\gppoint{gp mark 12}{(6.375,1.788)}
\gppoint{gp mark 12}{(6.378,1.788)}
\gppoint{gp mark 12}{(6.385,2.516)}
\gppoint{gp mark 12}{(6.388,1.788)}
\gppoint{gp mark 12}{(6.391,1.788)}
\gppoint{gp mark 12}{(6.394,1.788)}
\gppoint{gp mark 12}{(6.401,2.320)}
\gppoint{gp mark 12}{(6.404,2.567)}
\gppoint{gp mark 12}{(6.407,1.788)}
\gppoint{gp mark 12}{(6.410,2.614)}
\gppoint{gp mark 12}{(6.413,2.320)}
\gppoint{gp mark 12}{(6.416,2.124)}
\gppoint{gp mark 12}{(6.420,1.985)}
\gppoint{gp mark 12}{(6.423,2.459)}
\gppoint{gp mark 12}{(6.426,2.878)}
\gppoint{gp mark 12}{(6.429,2.395)}
\gppoint{gp mark 12}{(6.432,2.730)}
\gppoint{gp mark 12}{(6.435,3.011)}
\gppoint{gp mark 12}{(6.439,2.694)}
\gppoint{gp mark 12}{(6.442,2.516)}
\gppoint{gp mark 12}{(6.445,3.145)}
\gppoint{gp mark 12}{(6.448,2.459)}
\gppoint{gp mark 12}{(6.451,1.453)}
\gppoint{gp mark 12}{(6.454,2.764)}
\gppoint{gp mark 12}{(6.461,1.985)}
\gppoint{gp mark 12}{(6.470,1.453)}
\gppoint{gp mark 12}{(6.477,1.453)}
\gppoint{gp mark 12}{(6.480,1.453)}
\gppoint{gp mark 12}{(6.483,1.453)}
\gppoint{gp mark 12}{(6.486,1.453)}
\gppoint{gp mark 12}{(6.489,1.453)}
\gppoint{gp mark 12}{(6.496,1.453)}
\gppoint{gp mark 12}{(6.505,1.453)}
\gppoint{gp mark 12}{(6.511,1.788)}
\gppoint{gp mark 12}{(6.515,1.453)}
\gppoint{gp mark 12}{(6.524,1.788)}
\gppoint{gp mark 12}{(6.527,1.788)}
\gppoint{gp mark 12}{(6.534,1.788)}
\gppoint{gp mark 12}{(6.537,1.788)}
\gppoint{gp mark 12}{(6.540,1.453)}
\gppoint{gp mark 12}{(6.546,1.985)}
\gppoint{gp mark 12}{(6.556,1.453)}
\gppoint{gp mark 12}{(6.562,1.453)}
\gppoint{gp mark 12}{(6.572,1.453)}
\gppoint{gp mark 12}{(6.575,1.788)}
\gppoint{gp mark 12}{(6.578,1.453)}
\gppoint{gp mark 12}{(6.581,1.453)}
\gppoint{gp mark 12}{(6.587,1.453)}
\gppoint{gp mark 12}{(6.591,1.453)}
\gppoint{gp mark 12}{(6.594,1.985)}
\gppoint{gp mark 12}{(6.603,1.453)}
\gppoint{gp mark 12}{(6.610,1.453)}
\gppoint{gp mark 12}{(6.613,1.453)}
\gppoint{gp mark 12}{(6.616,1.788)}
\gppoint{gp mark 12}{(6.619,1.788)}
\gppoint{gp mark 12}{(6.622,1.788)}
\gppoint{gp mark 12}{(6.629,1.453)}
\gppoint{gp mark 12}{(6.632,1.453)}
\gppoint{gp mark 12}{(6.641,1.453)}
\gppoint{gp mark 12}{(6.645,1.453)}
\gppoint{gp mark 12}{(6.654,1.453)}
\gppoint{gp mark 12}{(6.664,1.453)}
\gppoint{gp mark 12}{(6.667,1.453)}
\gppoint{gp mark 12}{(6.670,1.453)}
\gppoint{gp mark 12}{(6.676,1.453)}
\gppoint{gp mark 12}{(6.683,1.453)}
\gppoint{gp mark 12}{(6.686,1.788)}
\gppoint{gp mark 12}{(6.692,1.453)}
\gppoint{gp mark 12}{(6.705,1.453)}
\gppoint{gp mark 12}{(6.714,1.453)}
\gppoint{gp mark 12}{(6.721,1.453)}
\gppoint{gp mark 12}{(6.724,1.788)}
\gppoint{gp mark 12}{(6.727,1.453)}
\gppoint{gp mark 12}{(6.730,1.453)}
\gppoint{gp mark 12}{(6.733,1.453)}
\gppoint{gp mark 12}{(6.736,1.985)}
\gppoint{gp mark 12}{(6.740,1.453)}
\gppoint{gp mark 12}{(6.749,1.453)}
\gppoint{gp mark 12}{(6.752,1.453)}
\gppoint{gp mark 12}{(6.755,1.453)}
\gppoint{gp mark 12}{(6.762,1.453)}
\gppoint{gp mark 12}{(6.765,1.453)}
\gppoint{gp mark 12}{(6.771,1.453)}
\gppoint{gp mark 12}{(6.778,1.453)}
\gppoint{gp mark 12}{(6.793,1.453)}
\gppoint{gp mark 12}{(6.797,1.788)}
\gppoint{gp mark 12}{(6.803,1.453)}
\gppoint{gp mark 12}{(6.806,1.453)}
\gppoint{gp mark 12}{(6.809,1.453)}
\gppoint{gp mark 12}{(6.819,4.904)}
\gpcolor{rgb color={0.000,0.000,0.000}}
\gpdefrectangularnode{gp plot 1}{\pgfpoint{1.196cm}{0.616cm}}{\pgfpoint{7.447cm}{5.075cm}}
\end{tikzpicture}
}\subfloat[Gramian $G(i,j)$ density distribution]{\centering{} 
\begin{tikzpicture}[scale=1.][gnuplot]
\path (0.000,0.000) rectangle (8.000,6.000);
\gpfill{rgb color={1.000,1.000,1.000}} (1.196,0.616)--(7.446,0.616)--(7.446,5.074)--(1.196,5.074)--cycle;
\gpcolor{color=gp lt color border}
\gpsetlinetype{gp lt border}
\gpsetlinewidth{.500}
\draw[gp path] (1.196,0.616)--(1.196,5.074)--(7.446,5.074)--(7.446,0.616)--cycle;
\gpsetlinewidth{0.50}
\draw[gp path] (1.196,0.616)--(1.447,0.616);
\draw[gp path] (7.447,0.616)--(7.196,0.616);
\gpcolor{rgb color={0.000,0.000,0.000}}
\node[gp node right,font={\fontsize{8pt}{9.6pt}\selectfont}] at (1.012,0.616) {$10^{-5}$};
\gpcolor{color=gp lt color border}
\draw[gp path] (1.196,0.884)--(1.321,0.884);
\draw[gp path] (7.447,0.884)--(7.322,0.884);
\draw[gp path] (1.196,1.041)--(1.321,1.041);
\draw[gp path] (7.447,1.041)--(7.322,1.041);
\draw[gp path] (1.196,1.153)--(1.321,1.153);
\draw[gp path] (7.447,1.153)--(7.322,1.153);
\draw[gp path] (1.196,1.239)--(1.321,1.239);
\draw[gp path] (7.447,1.239)--(7.322,1.239);
\draw[gp path] (1.196,1.310)--(1.321,1.310);
\draw[gp path] (7.447,1.310)--(7.322,1.310);
\draw[gp path] (1.196,1.370)--(1.321,1.370);
\draw[gp path] (7.447,1.370)--(7.322,1.370);
\draw[gp path] (1.196,1.421)--(1.321,1.421);
\draw[gp path] (7.447,1.421)--(7.322,1.421);
\draw[gp path] (1.196,1.467)--(1.321,1.467);
\draw[gp path] (7.447,1.467)--(7.322,1.467);
\draw[gp path] (1.196,1.508)--(1.447,1.508);
\draw[gp path] (7.447,1.508)--(7.196,1.508);
\gpcolor{rgb color={0.000,0.000,0.000}}
\node[gp node right,font={\fontsize{8pt}{9.6pt}\selectfont}] at (1.012,1.508) {$10^{-4}$};
\gpcolor{color=gp lt color border}
\draw[gp path] (1.196,1.776)--(1.321,1.776);
\draw[gp path] (7.447,1.776)--(7.322,1.776);
\draw[gp path] (1.196,1.933)--(1.321,1.933);
\draw[gp path] (7.447,1.933)--(7.322,1.933);
\draw[gp path] (1.196,2.045)--(1.321,2.045);
\draw[gp path] (7.447,2.045)--(7.322,2.045);
\draw[gp path] (1.196,2.131)--(1.321,2.131);
\draw[gp path] (7.447,2.131)--(7.322,2.131);
\draw[gp path] (1.196,2.202)--(1.321,2.202);
\draw[gp path] (7.447,2.202)--(7.322,2.202);
\draw[gp path] (1.196,2.261)--(1.321,2.261);
\draw[gp path] (7.447,2.261)--(7.322,2.261);
\draw[gp path] (1.196,2.313)--(1.321,2.313);
\draw[gp path] (7.447,2.313)--(7.322,2.313);
\draw[gp path] (1.196,2.359)--(1.321,2.359);
\draw[gp path] (7.447,2.359)--(7.322,2.359);
\draw[gp path] (1.196,2.400)--(1.447,2.400);
\draw[gp path] (7.447,2.400)--(7.196,2.400);
\gpcolor{rgb color={0.000,0.000,0.000}}
\node[gp node right,font={\fontsize{8pt}{9.6pt}\selectfont}] at (1.012,2.400) {$10^{-3}$};
\gpcolor{color=gp lt color border}
\draw[gp path] (1.196,2.668)--(1.321,2.668);
\draw[gp path] (7.447,2.668)--(7.322,2.668);
\draw[gp path] (1.196,2.825)--(1.321,2.825);
\draw[gp path] (7.447,2.825)--(7.322,2.825);
\draw[gp path] (1.196,2.937)--(1.321,2.937);
\draw[gp path] (7.447,2.937)--(7.322,2.937);
\draw[gp path] (1.196,3.023)--(1.321,3.023);
\draw[gp path] (7.447,3.023)--(7.322,3.023);
\draw[gp path] (1.196,3.094)--(1.321,3.094);
\draw[gp path] (7.447,3.094)--(7.322,3.094);
\draw[gp path] (1.196,3.153)--(1.321,3.153);
\draw[gp path] (7.447,3.153)--(7.322,3.153);
\draw[gp path] (1.196,3.205)--(1.321,3.205);
\draw[gp path] (7.447,3.205)--(7.322,3.205);
\draw[gp path] (1.196,3.251)--(1.321,3.251);
\draw[gp path] (7.447,3.251)--(7.322,3.251);
\draw[gp path] (1.196,3.291)--(1.447,3.291);
\draw[gp path] (7.447,3.291)--(7.196,3.291);
\gpcolor{rgb color={0.000,0.000,0.000}}
\node[gp node right,font={\fontsize{8pt}{9.6pt}\selectfont}] at (1.012,3.291) {$10^{-2}$};
\gpcolor{color=gp lt color border}
\draw[gp path] (1.196,3.560)--(1.321,3.560);
\draw[gp path] (7.447,3.560)--(7.322,3.560);
\draw[gp path] (1.196,3.717)--(1.321,3.717);
\draw[gp path] (7.447,3.717)--(7.322,3.717);
\draw[gp path] (1.196,3.828)--(1.321,3.828);
\draw[gp path] (7.447,3.828)--(7.322,3.828);
\draw[gp path] (1.196,3.915)--(1.321,3.915);
\draw[gp path] (7.447,3.915)--(7.322,3.915);
\draw[gp path] (1.196,3.985)--(1.321,3.985);
\draw[gp path] (7.447,3.985)--(7.322,3.985);
\draw[gp path] (1.196,4.045)--(1.321,4.045);
\draw[gp path] (7.447,4.045)--(7.322,4.045);
\draw[gp path] (1.196,4.097)--(1.321,4.097);
\draw[gp path] (7.447,4.097)--(7.322,4.097);
\draw[gp path] (1.196,4.142)--(1.321,4.142);
\draw[gp path] (7.447,4.142)--(7.322,4.142);
\draw[gp path] (1.196,4.183)--(1.447,4.183);
\draw[gp path] (7.447,4.183)--(7.196,4.183);
\gpcolor{rgb color={0.000,0.000,0.000}}
\node[gp node right,font={\fontsize{8pt}{9.6pt}\selectfont}] at (1.012,4.183) {$10^{-1}$};
\gpcolor{color=gp lt color border}
\draw[gp path] (1.196,4.452)--(1.321,4.452);
\draw[gp path] (7.447,4.452)--(7.322,4.452);
\draw[gp path] (1.196,4.609)--(1.321,4.609);
\draw[gp path] (7.447,4.609)--(7.322,4.609);
\draw[gp path] (1.196,4.720)--(1.321,4.720);
\draw[gp path] (7.447,4.720)--(7.322,4.720);
\draw[gp path] (1.196,4.807)--(1.321,4.807);
\draw[gp path] (7.447,4.807)--(7.322,4.807);
\draw[gp path] (1.196,4.877)--(1.321,4.877);
\draw[gp path] (7.447,4.877)--(7.322,4.877);
\draw[gp path] (1.196,4.937)--(1.321,4.937);
\draw[gp path] (7.447,4.937)--(7.322,4.937);
\draw[gp path] (1.196,4.989)--(1.321,4.989);
\draw[gp path] (7.447,4.989)--(7.322,4.989);
\draw[gp path] (1.196,5.034)--(1.321,5.034);
\draw[gp path] (7.447,5.034)--(7.322,5.034);
\draw[gp path] (1.196,5.075)--(1.447,5.075);
\draw[gp path] (7.447,5.075)--(7.196,5.075);
\gpcolor{rgb color={0.000,0.000,0.000}}
\node[gp node right,font={\fontsize{8pt}{9.6pt}\selectfont}] at (1.012,5.075) {$10^{0}$};
\gpcolor{color=gp lt color border}
\draw[gp path] (1.196,0.616)--(1.196,0.867);
\draw[gp path] (1.196,5.075)--(1.196,4.824);
\gpcolor{rgb color={0.000,0.000,0.000}}
\node[gp node center,font={\fontsize{8pt}{9.6pt}\selectfont}] at (1.196,0.308) {0};
\gpcolor{color=gp lt color border}
\draw[gp path] (2.446,0.616)--(2.446,0.867);
\draw[gp path] (2.446,5.075)--(2.446,4.824);
\gpcolor{rgb color={0.000,0.000,0.000}}
\node[gp node center,font={\fontsize{8pt}{9.6pt}\selectfont}] at (2.446,0.308) {0.2};
\gpcolor{color=gp lt color border}
\draw[gp path] (3.696,0.616)--(3.696,0.867);
\draw[gp path] (3.696,5.075)--(3.696,4.824);
\gpcolor{rgb color={0.000,0.000,0.000}}
\node[gp node center,font={\fontsize{8pt}{9.6pt}\selectfont}] at (3.696,0.308) {0.4};
\gpcolor{color=gp lt color border}
\draw[gp path] (4.947,0.616)--(4.947,0.867);
\draw[gp path] (4.947,5.075)--(4.947,4.824);
\gpcolor{rgb color={0.000,0.000,0.000}}
\node[gp node center,font={\fontsize{8pt}{9.6pt}\selectfont}] at (4.947,0.308) {0.6};
\gpcolor{color=gp lt color border}
\draw[gp path] (6.197,0.616)--(6.197,0.867);
\draw[gp path] (6.197,5.075)--(6.197,4.824);
\gpcolor{rgb color={0.000,0.000,0.000}}
\node[gp node center,font={\fontsize{8pt}{9.6pt}\selectfont}] at (6.197,0.308) {0.8};
\gpcolor{color=gp lt color border}
\draw[gp path] (7.447,0.616)--(7.447,0.867);
\draw[gp path] (7.447,5.075)--(7.447,4.824);
\gpcolor{rgb color={0.000,0.000,0.000}}
\node[gp node center,font={\fontsize{8pt}{9.6pt}\selectfont}] at (7.447,0.308) {1};
\gpcolor{color=gp lt color border}
\draw[gp path] (1.196,5.075)--(1.196,0.616)--(7.447,0.616)--(7.447,5.075)--cycle;
\gpcolor{rgb color={0.70,0.75,.71}}
\gpsetlinewidth{2.500}
\gpsetpointsize{8}
\gppoint{gp mark 12}{(1.227,4.929)}
\gppoint{gp mark 12}{(1.290,3.086)}
\gppoint{gp mark 12}{(1.352,3.109)}
\gppoint{gp mark 12}{(1.415,2.993)}
\gppoint{gp mark 12}{(1.477,3.032)}
\gppoint{gp mark 12}{(1.540,3.124)}
\gppoint{gp mark 12}{(1.602,3.019)}
\gppoint{gp mark 12}{(1.665,3.218)}
\gppoint{gp mark 12}{(1.727,3.178)}
\gppoint{gp mark 12}{(1.790,3.251)}
\gppoint{gp mark 12}{(1.852,3.858)}
\gppoint{gp mark 12}{(1.915,3.758)}
\gppoint{gp mark 12}{(1.977,3.600)}
\gppoint{gp mark 12}{(2.040,3.538)}
\gppoint{gp mark 12}{(2.102,3.450)}
\gppoint{gp mark 12}{(2.165,3.459)}
\gppoint{gp mark 12}{(2.227,3.373)}
\gppoint{gp mark 12}{(2.290,3.233)}
\gppoint{gp mark 12}{(2.352,3.170)}
\gppoint{gp mark 12}{(2.415,2.968)}
\gppoint{gp mark 12}{(2.477,3.117)}
\gppoint{gp mark 12}{(2.540,3.197)}
\gppoint{gp mark 12}{(2.602,2.860)}
\gppoint{gp mark 12}{(2.665,2.820)}
\gppoint{gp mark 12}{(2.727,2.856)}
\gppoint{gp mark 12}{(2.790,2.769)}
\gppoint{gp mark 12}{(2.853,2.820)}
\gppoint{gp mark 12}{(2.915,2.597)}
\gppoint{gp mark 12}{(2.978,2.513)}
\gppoint{gp mark 12}{(3.040,2.641)}
\gppoint{gp mark 12}{(3.103,2.389)}
\gppoint{gp mark 12}{(3.165,2.769)}
\gppoint{gp mark 12}{(3.228,2.893)}
\gppoint{gp mark 12}{(3.290,2.703)}
\gppoint{gp mark 12}{(3.353,2.405)}
\gppoint{gp mark 12}{(3.415,2.757)}
\gppoint{gp mark 12}{(3.478,2.587)}
\gppoint{gp mark 12}{(3.540,2.476)}
\gppoint{gp mark 12}{(3.603,2.232)}
\gppoint{gp mark 12}{(3.665,2.389)}
\gppoint{gp mark 12}{(3.728,2.121)}
\gppoint{gp mark 12}{(3.790,1.159)}
\gppoint{gp mark 12}{(3.853,1.852)}
\gppoint{gp mark 12}{(3.915,1.964)}
\gppoint{gp mark 12}{(3.978,1.695)}
\gppoint{gp mark 12}{(4.040,2.121)}
\gppoint{gp mark 12}{(4.103,1.159)}
\gppoint{gp mark 12}{(4.228,1.852)}
\gppoint{gp mark 12}{(4.290,2.319)}
\gppoint{gp mark 12}{(4.415,2.121)}
\gppoint{gp mark 12}{(4.478,2.121)}
\gppoint{gp mark 12}{(4.540,2.181)}
\gppoint{gp mark 12}{(4.603,2.278)}
\gppoint{gp mark 12}{(4.665,2.152)}
\gppoint{gp mark 12}{(4.728,2.232)}
\gppoint{gp mark 12}{(4.790,1.695)}
\gppoint{gp mark 12}{(4.853,2.050)}
\gppoint{gp mark 12}{(4.915,2.050)}
\gppoint{gp mark 12}{(5.103,1.695)}
\gppoint{gp mark 12}{(5.228,1.852)}
\gppoint{gp mark 12}{(5.290,1.964)}
\gppoint{gp mark 12}{(5.353,1.964)}
\gppoint{gp mark 12}{(5.415,1.852)}
\gppoint{gp mark 12}{(5.478,1.695)}
\gppoint{gp mark 12}{(5.540,2.050)}
\gppoint{gp mark 12}{(5.728,1.695)}
\gppoint{gp mark 12}{(5.790,2.050)}
\gppoint{gp mark 12}{(5.853,2.121)}
\gppoint{gp mark 12}{(5.978,1.852)}
\gppoint{gp mark 12}{(6.041,1.964)}
\gppoint{gp mark 12}{(6.103,1.852)}
\gppoint{gp mark 12}{(6.166,1.964)}
\gppoint{gp mark 12}{(6.228,2.232)}
\gppoint{gp mark 12}{(6.291,2.121)}
\gppoint{gp mark 12}{(6.353,2.389)}
\gppoint{gp mark 12}{(6.478,2.050)}
\gppoint{gp mark 12}{(6.541,1.852)}
\gppoint{gp mark 12}{(6.603,1.964)}
\gppoint{gp mark 12}{(6.666,1.695)}
\gppoint{gp mark 12}{(6.728,1.695)}
\gppoint{gp mark 12}{(6.791,1.695)}
\gppoint{gp mark 12}{(6.853,2.050)}
\gppoint{gp mark 12}{(6.978,1.964)}
\gppoint{gp mark 12}{(7.041,1.695)}
\gppoint{gp mark 12}{(7.103,2.389)}
\gppoint{gp mark 12}{(7.291,1.964)}
\gppoint{gp mark 12}{(7.353,2.181)}
\gppoint{gp mark 12}{(7.416,3.148)}
\gpcolor{rgb color={0.000,0.000,0.000}}
\gpdefrectangularnode{gp plot 1}{\pgfpoint{1.196cm}{0.616cm}}{\pgfpoint{7.447cm}{5.075cm}}
\end{tikzpicture}
}%
}\protect\caption{\label{fig:fanbeam-1}(a) Example normalized angles distribution (y-axis)
for the first eight columns of $AA^{T}$ where $\theta_{i,j}=\cos^{-1}(\langle\hat{\bm{a}}_{i},\hat{\bm{a}}_{j}\rangle)\,\,\forall\,\, i,j\,\in\left\{ 1,\ldots,M\right\} $
vs angles (x-axis) degrees using fan-beam tomographic data acquisition
strategy, (b) Gramian matrix $\langle\hat{\bm{a}}_{i},\hat{\bm{a}}_{j}\rangle$
distribution}
\end{figure}

\subsection{Angular Distribution of Hyperplanes}

A comparison of the distribution of hyperplane sampling angles in
computed tomography (CT) was performed to investigate the convergence
rate versus measurement strategy. Example results are presented for
iterative convergence of methods K, RK, and RKHA under conditions
of random, fan, and parallel beam sampling strategies using the Shepp-Logan
phantom (see Figure \eqref{Sta:Phantom-image-1})%
\footnote{Shepp-Logan phantom was generated from AIRtools/paralleltomo.m with
non-uniform coherent parallel tomographic CT sampling, P. C. Hansen
and M. Saxild-Hansen, AIR Tools - A MATLAB (tm) Package of Algebraic
Iterative Reconstruction Methods, Journal of Computational and Applied
Mathematics, 236 (2012), pp. 2167-2178%
}, \emph{paralleltomo.m} and \emph{fanbeamtomo.m} from the AIRtools
distribution \cite{Hansen20122167}, and \emph{randn()} from the built-in
function method \cite{Marsaglia:Tsang:2000:JSSOBK:v05i08}.

\subsection{Measurement Coherence}

In linear algebra, the coherence or mutual coherence \cite{1564423}
of a row measurement matrix $A$ is defined as the maximum absolute
value of the cross-correlations between the normalized rows of $A$.

Formally, let $\left\{ \bm{a}_{1},\ldots,\bm{a}_{M}\right\} \in{\mathbb{\R}}^{N}$
be the set of row vectors of the matrix $A\in\R^{M\times N}$  normalized
such that $\langle\bm{a}_{i},\bm{a}_{i}\rangle=\bm{a}_{i}^{H}\bm{a}_{i}=1$
where $(.)^{H}$ is the Hermitian conjugate and where $M>N.$ Let
the mutual coherence of $A$ be defined as 
\begin{equation}
\phi_{i,j}=\max_{1\le i\ne j\le M}\left|\bm{a}_{i}^{H}\bm{a}_{j}\right|.\label{eq:coherence}
\end{equation}
A lower bound was derived as $\phi\ge\frac{M-N}{N(M-1)}$ in reference
Welch \cite{1055219}.

It is noted that the statistical expectation%
\footnote{A more formal treatment of the expectation of random IID vectors is
given in section \eqref{sec:Consider-the-expectationofIID}.%
} of the non-diagonal Gramian matrix elements $G_{i,j}=\langle\hat{\bm{a}}_{i},\hat{\bm{a}}_{j}\rangle\,(1\le i\ne j\le M)$
for normalized random unit vectors $\left\{ \hat{\bm{a}}_{i},\hat{\bm{a}}_{j}\right\} $
would be zero for two independent random IID row vectors, $1/N$
for the case of a single dependent vector component (one variable
in $N$ variables), and the maximum expected value occurs when two
unit row vectors are parallel, which gives a value of unity. Estimated
numerical results for the three sampling methods are shown in Table
\eqref{tab:Computed-coherence-for} along with values for the mean
of the Gramian.
\begin{table}
\noindent \begin{centering}
\begin{tabular}{|c|c|c|c|}
\hline 
Coherence vs Measurement Method%
\footnote{sampling over range of $[0,2\pi]$ radians%
}  & Random & Fan & Parallel\tabularnewline
\hline 
\hline 
coherence Eq. \eqref{eq:coherence} & .4 & 1.0 & 1.0\tabularnewline
\hline 
average value of $G_{i,j}=\langle\hat{\bm{a}}_{i},\hat{\bm{a}}_{j}\rangle\,(1\le i\ne j\le M)$ & -.0013 & .06 & .18\tabularnewline
\hline 
median value of $G_{i,j}=\langle\hat{\bm{a}}_{i},\hat{\bm{a}}_{j}\rangle\,(1\le i\ne j\le M)$ & -0.0009 & 0 & .12\tabularnewline
\hline 
\end{tabular}
\par\end{centering}

\protect\caption{\label{tab:Computed-coherence-for}Typical coherence estimates for
$N=100,\,\, M=200$ for random \emph{randn()} and $N=100,\,\, M=222$
for fan \emph{fanbeamtomo()} and parallel \emph{paralleltomo()}}
\end{table}

\begin{figure}
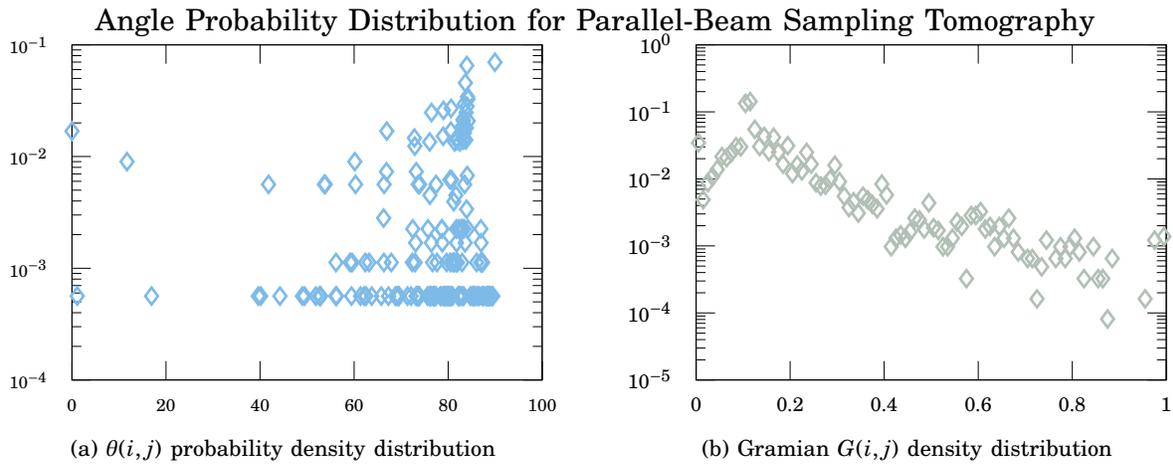

\noindent \centering{}Angle Probability Distribution for Parallel-Beam
Sampling Tomography\vspace{-33pt}
\makebox[1\columnwidth]{%
\subfloat[$\theta(i,j)$ probability density distribution ]{\centering{} % ctangles-compare-1-1.tex

}%
}\protect\caption{\label{fig:parallel}(a) Example normalized angles distribution (y-axis)
from $AA^{T}$ where $\theta_{i,j}=\cos^{-1}(\langle\hat{\bm{a}}_{i},\hat{\bm{a}}_{j}\rangle)\,\,\forall\,\, i,j\,\in\{1,\ldots,M\}$
vs angles (x-axis) degrees using parallel-beam tomographic data acquisition
strategy, (b) Gramian matrix $\langle\hat{\bm{a}}_{i},\hat{\bm{a}}_{j}\rangle$
distribution}
\end{figure}

Computations of the Gramian and angular density distributions are
shown in Figures \eqref{fig:random-1-2-1-1}, \eqref{fig:fanbeam-1},
and \eqref{fig:parallel}. It should be noted that the random sampling
is concentrated near $90$ degrees probability and zero for the Gramian,
but parallel sampling is spread out across the interval $[0,90]$
degrees.

\subsection{Distribution of Measurement Angles for K, RK, and RKHA for Shepp-Logan
Versus Measurement Method}

Firstly, the convergence rates of K, RK, and RKHA are noted to be
closely correlated for the case of random data sampling of the phantom.
This is consistent with the mean values of coherence near zero for
random sampling.

The cases for fan and parallel sampling have increasingly higher coherence,
and generally benefit from methods which minimize the coherence, such
as RK, RKHA, and RKOS.

Representative results for convergence are shown in Figures \eqref{fig:Example-convergence-result random},
\eqref{fig:Example-convergence-result fan}, and \eqref{fig:Example-convergence-result para}.
Comparison of convergence results to the estimated coherence for the
three cases given in Table \eqref{tab:Computed-coherence-for} suggest
consistent interpretation.

Since the iterative methods utilize projections, the angles between
the optical lines of sight (LOS) forming the measurement hyperplanes
is of considerable interest. The figures also show example computations
of distribution of measurement hyperplane angles relative to a hyperplane
reference as given by $\theta_{i,j}=\cos^{-1}(\langle\hat{\bm{a}}_{i},\hat{\bm{a}}_{j}\rangle)\,\,\forall\,\,\hat{\bm{a}}_{i},\hat{\bm{a}}_{j}\,\in A_{i}\,\,\forall i,j\in\left\{ 1,\ldots,M\right\} $
where the unit norm vectors $\hat{\bm{a}}_{i},\hat{\bm{a}}_{j}$ are
selected rows of $A$.

\subsection{Convergence of K, RK, and RKHA vs Measurement Method}

\begin{figure}
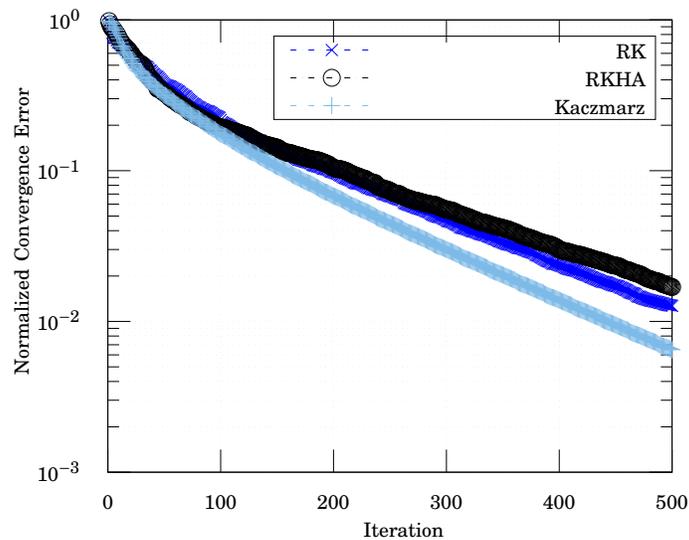

\centering{} 
%
\protect\caption{\label{fig:Example-convergence-result random}Semilog (y-axis) plot
example convergence result for K, RK, RKHA on Shepp-Logan phantom
using IID random tomographic data acquisition for 10 cycles of iteration
(x-axis). Note that randomization tends to equalize convergence rates
and diminish advantage of a particular method.}
\end{figure}
\begin{figure}
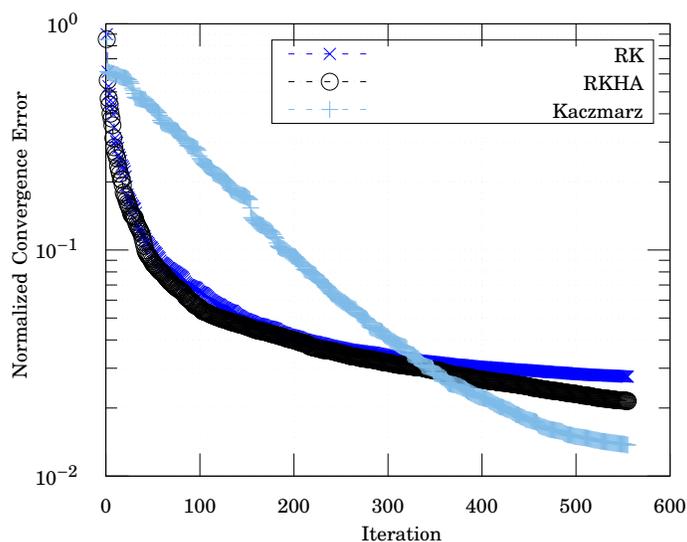

\centering{} 

\protect\caption{\label{fig:Example-convergence-result para}Semilog (y-axis) plot
example convergence result for K, RK, RKHA on Shepp-Logan phantom
using parallel tomographic data acquisition for 10 cycles of iteration
(x-axis). Note that initially, both RK and RKHA have similar advantage,
but simple Kaczmarz eventually improves.}
\end{figure}

Iterative simulations were performed to estimate the relative convergence
rates of methods K, RK, RKHA for the data examples above, random,
parallel, and fan beam sampling.
\begin{figure}
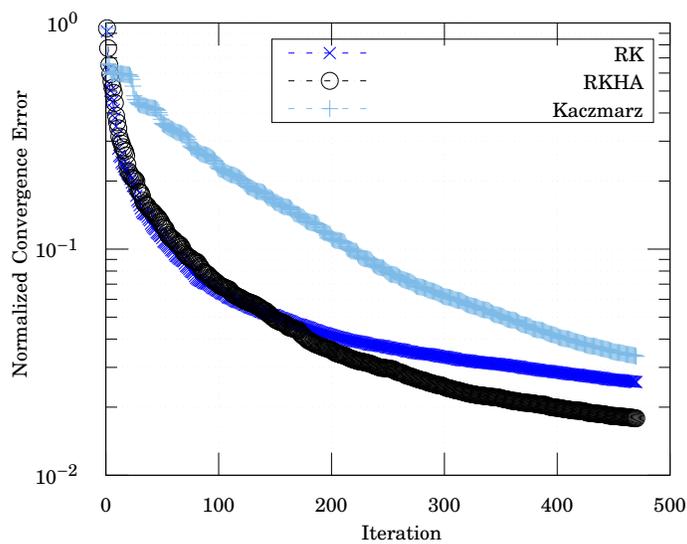

\centering{} 
%
\protect\caption{\label{fig:Example-convergence-result fan}Semilog (y-axis) plot example
convergence result for K, RK, RKHA on Shepp-Logan phantom using fan
tomographic data acquisition for 10 cycles of iteration (x-axis).
Note that both RK and RKHA appear to have advantage since each method
utilizes randomization which improves avoidance of coherent neighbors,
but simple Kaczmarz is too naive.}
\end{figure}
 Representative results are shown in Figures \eqref{fig:Example-convergence-result random},
\eqref{fig:Example-convergence-result fan}, and \eqref{fig:Example-convergence-result para}
for noiseless data measurement scenarios of the standard Shepp-Logan
phantom. 

\section{Conclusions}

A new iterative selection rule based upon the relative central angle
shows enhanced convergence in measurements which contain coherence.
However, the method requires a computational penalty related to the
dot products of all to all rows, which may be overcome by \emph{a
priori} determination. A new block method using constructed orthogonal
subspace projections provides enhanced tolerance to measurement incoherence,
but may be affected by noise at least as much as simple Kaczmarz.
The exponential convergence is accelerated by the $P/N$ term and
is computationally feasible for small $P$ relative to $N$.

The convergence of above subspace methods was demonstrated using statistical
IID assumptions. But, the more generalized approach based upon cyclical
projections using the formalism of Galantai also prove convergence,
without the statistical argument. 

It is worthwhile to note that an additional method to prove the convergence
rate for a given angular probability distribution function is currently
underway and is considered an essential task towards validation of
the RKHA results.

\section*{Acknowledgments}

The authors thank Akram Aldroubi and Alex Powell for their invaluable
feedback. The research of Ali Sekmen is supported in part by NASA
Grant NNX12AI14A.

\bibliographystyle{IEEEtran}

\end{document}